\documentclass[aap,MSNbibl,citesort,dvips]{arximspdf}
\usepackage{mathbh}
\usepackage{graphicx}

\doi{10.1214/11-AAP803}
\volume{22}
\issue{4}
\pubyear{2012}
\firstpage{1611}
\lastpage{1641}

\makeatletter

\newtheorem{lemma}{Lemma}[section]
\newtheorem{theorem}{Theorem}[section]

\newcommand{\ttilde}[1]{\hspace*{1pt}\tilde{\hspace*{-1pt}\tilde{#1}}{}}

\newcommand{\N}{{\mathbb{N}}}
\newcommand{\R}{{\mathbb{R}}}
\newcommand{\E}{{\mathbb{E}}}
\newcommand{\ld}{{\lambda}}

\newcommand{\norm}[1]{{\Arrowvert #1\Arrowvert}}

\makeatother

\begin{document}
\begin{frontmatter}

\title{Strong convergence of an explicit numerical method for
SDEs with nonglobally Lipschitz continuous coefficients\thanksref{T1}}
\runtitle{SDEs with nonglobally Lipschitz continuous coefficients}

\thankstext{T1}{Supported in part by the Collaborative
Research Centre 701 ``Spectral Structures and Topological Methods
in Mathematics'' and by the research project ``Numerical solutions
of stochastic differential equations with nonglobally Lipschitz
continuous coefficients'' both funded by the
German Research Foundation.}

\begin{aug}
\author[A]{\fnms{Martin}~\snm{Hutzenthaler}\ead[label=e1]{hutzenthaler@bio.lmu.de}},
\author[B]{\fnms{Arnulf}~\snm{Jentzen}\ead[label=e2]{ajentzen@math.princeton.edu}}
\and
\author[C]{\fnms{Peter~E.}~\snm{Kloeden}\corref{}\ead[label=e3]{kloeden@math.uni-frankfurt.de}}
\runauthor{M. Hutzenthaler, A. Jentzen
and P. E. Kloeden}
\affiliation{University of Munich, Princeton University and Goethe University}
\address[A]{M. Hutzenthaler\\
LMU Biozentrum\\
Department Biologie II\\
University of Munich (LMU)\\
D-82152 Planegg-Martinsried\\
Germany\\
\printead{e1}} 
\address[B]{A. Jentzen\\
Program in Applied\\
\quad and Computational Mathematics\\
Princeton University\\
Princeton, New Jersey 08544-1000\\
USA\\
\printead{e2}}
\address[C]{P. E. Kloeden\\
Institute for Mathematics\\
Goethe University Frankfurt am Main\\
D-60054, Frankfurt am Main\\
Germany\\
\printead{e3}}
\end{aug}

\received{\smonth{3} \syear{2011}}

%
\begin{abstract}
On the one hand, the explicit Euler scheme fails to converge strongly
to the exact solution of a stochastic differential equation (SDE) with
a superlinearly growing and globally one-sided Lipschitz continuous
drift coefficient. On the other hand, the implicit Euler scheme is
known to converge strongly to the exact solution of such an SDE.
Implementations of the implicit Euler scheme, however, require
additional computational effort. In this article we therefore propose
an explicit and easily implementable numerical method for such an SDE
and show that this method converges strongly with the standard order
one-half to the exact solution of the SDE. Simulations reveal that this
explicit strongly convergent numerical scheme is considerably faster
than the implicit Euler scheme.
\end{abstract}

%

\end{frontmatter}

\section{Introduction and main result}\label{secintro}
The explicit
Euler scheme
(see, e.g., Kloeden and Platen~\cite{kp92}, Maruyama~\cite{m55}
and Milstein~\cite{m95})
is most commonly used for approximating stochastic differential
equations (SDEs) with globally Lipschitz continuous coefficients.
Unfortunately, the explicit Euler scheme does not converge in the
strong mean square sense to the exact solution of an SDE with a
superlinearly growing and globally one-sided Lipschitz continuous drift
coefficient. Even worse, Theorem 1 in~\cite{hjk11} shows for such an
SDE that the absolute moments of the explicit Euler approximations at a
finite time point $T\in(0,\infty)$\vadjust{\goodbreak} diverge to infinity. The implicit
Euler scheme is better than the explicit Euler scheme in that it
converges strongly to the exact solution of such an SDE (see Higham,
Mao and Stuart~\cite{hms02}). However, additional computational effort
is required for its implementation. Therefore, we wish to identify
explicit numerical methods which are strongly convergent even for SDEs
with superlinearly growing coefficients. In this article we propose a
``tamed'' version of the explicit Euler scheme in which the drift term
is modified such that it is uniformly bounded. Being almost identical
to the explicit Euler method, this version is explicit and easy to
implement. Now the benefit of this ``tamed'' Euler scheme is that it
converges strongly to the exact solution in case of SDEs with
superlinearly growing coefficients. More precisely, the main result of
this article shows that this ``tamed'' Euler scheme converges strongly
with the standard convergence order $ \frac{1}{2} $ to the exact
solution of the SDE if the drift coefficient function is globally
one-sided Lipschitz continuous and has an at most polynomially growing
derivative. The diffusion coefficient is assumed to be globally
Lipschitz continuous here. Simulations confirm our theoretical results.

Throughout the whole
article we assume that the following setting is fulfilled.
Let $ T \in(0,\infty) $ be a fixed real number, let $ ( \Omega,
\mathcal{F}, \mathbb{P} ) $ be a probability space\vspace*{1pt} with normal
filtration $ ( \mathcal{F}_t )_{ t \in[0,T] } $, let $ d, m
\in\mathbb{N} := \{1,2,\ldots\} $, let $ W = ( W^{ (1) }, \ldots, W^{ (m)
} ) \dvtx[0,T] \times\Omega\rightarrow\mathbb{R}^m $ be an $ m
$-dimensional standard $ ( \mathcal{F}_t )_{ t \in[0,T] } $-Brownian
motion and let $ \xi\dvtx\Omega\rightarrow\mathbb{R}^d $ be an $
\mathcal{F}_0 / \mathcal{B}( \mathbb{R}^d ) $-measurable mapping with
$ \mathbb{E}[ \| \xi\|^p ] < \infty$ for all $ p \in[1,\infty) $. Here
and below we use the notation $ \| v \| := ( | v_1 |^2 + \cdots+ | v_k
|^2 )^{{1/2}} $, $ \langle v, w \rangle := v_1 \cdot w_1 + \cdots+ v_k
\cdot w_k $ for all $ v = (v_1, \ldots, v_k) $, $ w = (w_1, \ldots, w_k)
\in\mathbb{R}^k $, $ k \in\mathbb{N} $, and $ \| A \| := {\sup_{ v
\in\mathbb{R}^l, \| v \| \leq1 }} \| A v \| $ for all $ A
\in\mathbb{R}^{ k \times l } $, $ k, l \in\mathbb{N} $. Moreover, let $
\mu\dvtx\mathbb{R}^d \rightarrow\mathbb{R}^d $ be a continuously
differentiable and globally one-sided Lipschitz continuous function
whose derivative grows at most polynomially and let $ \sigma= (
\sigma_{i,j} )_{
i \in\{ 1, 2, \ldots, d \}, j \in\{ 1, 2, \ldots, m \} }
\dvtx\mathbb{R}^d \rightarrow\mathbb{R}^{ d \times m } $ be a globally
Lipschitz continuous function. More formally, suppose that there is a
real number $ c \in(0,\infty) $ such that $ \| \mu'( x ) \| \leq c( 1 +
\| x \|^c ) $, $ \| \sigma( x ) - \sigma( y ) \| \leq c \| x - y \| $
and $ \langle x-y, \mu( x ) - \mu( y ) \rangle \leq c \| x - y \|^2 $ for all $ x,
y \in\mathbb{R}^d $. Then consider the SDE
%
%
\begin{equation}
\label{eqSDE} dX_t = \mu( X_t ) \,dt + \sigma( X_t ) \,dW_t,\qquad X_0 =
\xi,
\end{equation}
for $ t \in[0,T] $. The drift coefficient $\mu$ is the infinitesimal
mean of the process~$X$ and the diffusion coefficient $\sigma$ is the
infinitesimal standard deviation of the process~$X$. Under the above
assumptions, the SDE~(\ref{eqSDE}) is known to have a~unique strong
solution. More formally, there exists an adapted stochastic process $ X
\dvtx[0,T] \times\Omega\rightarrow\mathbb{R}^d $ with continuous
sample paths fulfilling
%
%
\begin{equation}
\label{eqSDEsol} X_t = \xi+ \int_0^t \mu( X_s ) \,ds + \int_0^t
\sigma( X_s ) \,dW_s
\end{equation}
for all $ t \in[0,T] $
$ \mathbb{P} $-a.s.
We refer to
Theorem 2 in Alyushina~\cite{Alyushina1987}, Theorem 1 in Krylov
\cite{Krylov1990} and Theorem 2.4.1 in Mao~\cite{m97} for existence and
uniqueness results for SDEs of the
form~(\ref{eqSDE}).\vadjust{\goodbreak} 

The goal of this article is to solve the strong approximation problem
(see, e.g., Kloeden and Platen~\cite{kp92}, Section 9.3) of the SDE
(\ref{eqSDE}). More precisely, our aim is to find a numerical
approximation $ Y \dvtx\Omega\rightarrow\mathbb{R}^d $ which satisfies
%
%
\begin{equation}
\label{eqstrongproblem} ( \mathbb{E} [ \| X_T - Y \|^2 ] )^{
{1/2} } < \varepsilon
\end{equation}
for a given precision $ \varepsilon> 0 $ and which can be implemented
with as little computational effort as possible. At this point let us
comment on the importance of solving the strong approximation problem
(\ref{eqstrongproblem}). A central motivation for studying strong
approximations in the sense of~(\ref{eqstrongproblem}) is Giles'
seminal paper~\cite{g08b} (see also Heinrich~\cite{h98}). There he
introduces, in comparison to the classical Monte Carlo method, a very
efficient, somehow accelerated Monte Carlo method for approximating
moments or other expectations of functionals of the SDE solution via
numerical schemes that converge strongly (see also Creutzig, Dereich,
M\"{u}ller-Gronbach and Ritter~\cite{cdmr09} for an detailed comparison
of the classical and the new ``accelerated'' Monte Carlo method). In
view of this method, strong approximations of the exact solution of the
SDE~(\ref{eqSDE}) in the sense of~(\ref{eqstrongproblem}) yield very
efficient approximations of expectations of functionals of the SDE
solution and this is a central reason for developing strongly
convergent numerical methods.

The simplest and most obvious idea to solve the strong approximation
problem~(\ref{eqstrongproblem}) is to apply the explicit Euler scheme
to the SDE~(\ref{eqSDE}). More precisely, the explicit Euler method for
the SDE~(\ref{eqSDE}) is given by mappings $ \tilde{Y}_n^N \dvtx\Omega
\rightarrow\mathbb{R}^d $, $ n \in\{ 0, 1, \ldots, N \} $, $ N
\in\mathbb{N}$, which satisfy $ \tilde{Y}_0^N = \xi$ and
%
%
\begin{equation}
\label{eqEulerscheme} \tilde{Y}_{ n + 1 }^N = \tilde{Y}_n^N + \frac{ T
}{ N } \cdot\mu( \tilde{Y}_n^N ) + \sigma( \tilde{Y}_n^N ) \bigl( W_{
{ (n + 1) T }/{ N } } - W_{ { n T }/{ N } } \bigr)
\end{equation}
for all $ n \in\{ 0, 1, \ldots, N - 1 \} $ and all $ N \in\mathbb{N}
$. In the literature (see, e.g., Theorem~10.2.2 in Kloeden and
Platen~\cite{kp92}, Theorem 1.1 in Milstein~\cite{m95} or Theorem~3.1
in Yuan and Mao~\cite{ym08}) the convergence results for the explicit
Euler scheme require the drift coefficient $\mu$ of the SDE
(\ref{eqSDE}) to be globally Lipschitz continuous or to grow at most
linearly, which we have not assumed in our setting.
As it turns out, the assumption of an at most linearly growing drift
function is essentially necessary. More precisely, in the case $ d=m=1
$, it has recently been shown in~\cite{hjk11} that the root mean square
distance of the exact solution of the SDE~(\ref{eqSDE}) and of the
explicit Euler approximation~(\ref{eqEulerscheme}) diverges to infinity
%
%
\begin{equation}
\label{eqEulerdiv} \lim_{ N \rightarrow\infty} ( \mathbb{E} [ \| X_T -
\tilde{Y}_N^N \|^2 ] )^{ {1/2} } = \infty,
\end{equation}
if the drift coefficient $ \mu$ of the SDE~(\ref{eqSDE}) grows
superlinearly, that is, if there are real numbers $ \alpha, C
\in(1,\infty) $ such that $ | \mu( x ) | \geq\frac{ | x |^{ \alpha} }{
C } $ holds for all $ | x | \geq C $. Thus the explicit Euler scheme
(\ref{eqEulerscheme}) does not solve the strong approximation problem
(\ref{eqstrongproblem}) of the SDE~(\ref{eqSDE}) in general. This is
particularly unfortunate as SDEs with superlinearly growing
coefficients are quite important in applications\vadjust{\goodbreak} (see, e.g., \cite
{BenguriaKac1981,GinzburgLandau1950,HutzenthalerWakolbinger2007,KhasminskiiKlebaner2001,Lythe1995,Oettinger1996,Schuss1980}).
We remark that in contrast to strong mean square
convergence, pathwise convergence of
the explicit Euler
method~(\ref{eqEulerscheme}) to the
exact solution of the SDE~(\ref{eqSDE}) holds
due to Gy{\"o}ngy's
result~\cite{g98b}.

Another idea for solving the strong approximation problem
(\ref{eqstrongproblem}) is to apply the implicit Euler scheme, a.k.a.
backward Euler scheme (see Higham, Mao and Stuart~\cite{hms02}), to the
SDE~(\ref{eqSDE}). The implicit Euler scheme for the SDE~(\ref{eqSDE})
is given by mappings $ \ttilde{{Y}}_n^N
\dvtx\Omega\rightarrow\mathbb{R}^d $, $ n \in\{ 0, 1, \ldots, N \} $, $
N \in\mathbb{N} $, which satisfy $ \ttilde{{Y}}_0^N = \xi$ and
%
%
\begin{equation}
\label{eqimplicitEulerscheme} \ttilde{{Y}}_{ n + 1 }^N =
\ttilde{{Y}}_n^N + \frac{ T }{ N } \cdot\mu( \ttilde{{Y}}_{
n + 1 }^N ) + \sigma( \ttilde{{Y}}_n^N ) \bigl( W_{ { (n + 1) T }/{
N } } - W_{ { n T }/{ N } } \bigr)
\end{equation}
for all $ n \in\{ 0, 1, \ldots, N - 1 \} $ and all $ N \in\mathbb
{N} $.
A solution of this implicit equation is guaranteed to exist and to be
unique for $N\in\mathbb{N}$ large enough due to the globally one-sided
Lipschitz continuity of $\mu$.
In the same setting as in this article, Higham, Mao and Stuart showed
in Theorem 5.3 in~\cite{hms02} (see also
\cite{h10,hmps10,h96,lae08,s10,sm10a,sm10b,t02} and the references
therein for more approximation results on implicit numerical methods
for SDEs of the form~(\ref{eqSDE})) that the implicit Euler scheme
(\ref{eqimplicitEulerscheme}) converges with order $ \frac{1}{2} $ to
the exact solution of the SDE~(\ref{eqSDE}) in the root mean square
sense, that is, they established the existence of a real number $ C
\in[0,\infty) $ such that
%
%
\begin{equation}
( \mathbb{E} [ \| X_T - \ttilde{{Y}}_N^N \|^2 ] )^{ {1/2} }
\leq C \cdot N^{ -{1}/{2} }
\end{equation}
for all $ N \in\mathbb{N} $. However, additional computational effort
is required in order to implement~(\ref{eqimplicitEulerscheme}) since
the zero of a nonlinear equation has to be determined in each time step
in~(\ref{eqimplicitEulerscheme}).

To sum up, the explicit Euler scheme, on the one hand, is explicit and
easily implemented but does, in general, not converge strongly to the
exact solution of the SDE~(\ref{eqSDE}). The implicit Euler scheme, on
the other hand, converges strongly to the exact solution of the SDE
(\ref{eqSDE}) but additional computational effort is required for its
implementation.
Therefore, we aim at a simple
explicit numerical method
which
converges strongly to the exact
solution of the
SDE~(\ref{eqSDE}).


More formally, the
following numerical method for approximating
the solution of the SDE~(\ref{eqSDE})
is proposed here.
Let $ Y^N_n \dvtx\Omega\rightarrow\mathbb{R}^d $, $ n
\in\{0,1,\ldots,N\}$, $N \in\mathbb{N}$, be given by $ Y^N_0 = \xi$ and
%
%
\begin{equation}
\label{eqscheme} Y_{ n+ 1 }^N = Y_n^N + \frac{ {T}/{N} \cdot\mu(
Y_n^N ) }{ 1 + {T}/{N} \cdot\| \mu( Y_n^N ) \| } + \sigma( Y_n^N )
\bigl( W_{ { (n + 1) T }/{ N } } - W_{ { n T }/{ N } } \bigr)
\end{equation}
for all $ n \in\{ 0,1,\ldots,N-1 \}$ and all $ N \in\mathbb{N} $. We
refer to the numerical meth\-od~(\ref{eqscheme}) as a \textit{tamed Euler
scheme}. In this method the drift term $ \frac{T}{N} \cdot\mu(Y_n^N)$
is ``tamed'' by the factor $ 1 / ( 1 + \frac{T}{N} \cdot\| \mu( Y_n^N
) \| ) $ for $ n \in\{ 0, 1, \ldots, N - 1 \} $ and $ N \in\mathbb{N} $
in~(\ref{eqscheme}). Note that the norm of $ \frac{T}{N}
\cdot\mu(Y_n^N)
/ ( 1 + \frac{T}{N} \cdot\| \mu( Y_n^N ) \| ) $ is bounded by $1$
for every $ n \in\{ 0, 1, \ldots, N - 1 \} $ and every\vadjust{\goodbreak} $ N
\in\mathbb{N} $. This\vspace*{1pt} prevents the drift term from producing
extraordinary large values.
Additionally, the Taylor expansion of $ \frac{T}{N} \cdot\mu(x) / ( 1
+ \frac{T}{N} \cdot\| \mu( x ) \| ) $ in $1/N$ for fixed
$x\in\mathbb{R}$ is equal to the drift term $ \frac{T}{N} \cdot\mu(x) $
plus terms of order $O(\frac{1}{N^2})$. More formally, we see that
%
%
\begin{eqnarray}
\label{eqsecondorder}
Y_{ n+ 1 }^N &=& Y_n^N + \frac{T}{N} \cdot\mu(
Y_n^N ) + \sigma( Y_n^N ) \bigl( W_{ { (n + 1) T }/{ N } } - W_{ {
n T }/{ N } } \bigr)\nonumber\\[-8pt]\\[-8pt]
&&{} - \biggl(\frac{T}{N}\biggr)^{ 2} \frac{\mu(Y_n^N) \cdot
\|\mu(Y_n^N)\|} { 1 + {T}/{N} \cdot\| \mu( Y_n^N ) \| }\nonumber
\end{eqnarray}
for all $ n \in\{ 0, 1, \ldots, N-1 \} $ and all $ N \in\mathbb{N} $.
Thus the tamed Euler scheme~(\ref{eqscheme}) coincides with the
explicit Euler method~(\ref{eqEulerscheme}) up to terms of second
order.
Moreover, note that the tamed Euler scheme~(\ref{eqscheme}) can be
simulated easily and the drift function $\mu$ needs to be evaluated
only once in each iteration of~(\ref{eqscheme}).
More precisely,
having calculated
$v:=\frac{T}{N}
\cdot\mu(Y_n^N)$,
the drift term in~(\ref{eqscheme})
is then readily computed as
$\frac{v}{1+\|v\|}$.\vspace*{1pt}

In order to formulate our convergence theorem for the tamed Euler
meth\-od~(\ref{eqscheme}), we now introduce appropriate time continuous
interpolations\vspace*{1pt} of the time discrete numerical approximations
(\ref{eqscheme}).
More formally, let $ \bar{Y}^N \dvtx[0,T] \times\Omega\rightarrow
\mathbb{R}^d $, $ N \in\mathbb{N} $, be a sequence of stochastic
processes given by
%
%
\begin{equation}
\label{eqtimecont} \bar{Y}_{ t }^N = Y_n^N + \frac{ ( t - { n T }/{
N } ) \cdot\mu( Y_n^N ) }{ 1 + {T}/{N} \cdot\| \mu( Y_n^N ) \| }
+ \sigma( Y_n^N ) ( W_{ t } - W_{ { n T }/{ N } } )
\end{equation}
for all $ t \in[ \frac{ n T }{ N }, \frac{ ( n+ 1) T }{ N } ] $, $ n
\in\{ 0,1 , \ldots, N-1 \} $ and all $ N \in\mathbb{N} $. Note that $
\bar{Y}^N \dvtx\break [0$, $T] \times\Omega\rightarrow\mathbb{R}^d $ is an
adapted stochastic process with continuous sample paths for every $ N
\in\mathbb{N} $. We are now ready to formulate the main result of this
article.
%
%
\begin{theorem}[(Main result)]
\label{thmmain}
\label{thmmainresult}
Let the setting in this section be fulfilled. Then there exists a
family $ C_p \in[ 0, \infty) $, $ p \in[1,\infty) $, of real numbers
such that
%
%
\begin{equation}
\label{eqmainresult} \Bigl( \mathbb{E} \Bigl[ {\sup_{ t \in[0,T] }} \| X_t -
\bar{Y}_t^N \|^p \Bigr] \Bigr)^{ {1/p} } \leq C_p \cdot N^{ -{1}/{2} }
\end{equation}
for all $ N \in\mathbb{N} $ and all $ p \in[1,\infty) $. Here $ X
\dvtx[0,T] \times\Omega\rightarrow\mathbb{R}^d $ is the exact
solution of the SDE~(\ref{eqSDE}) and $ \bar{Y}^N \dvtx[0,T] \times
\Omega\rightarrow\mathbb{R}^d $, $ N \in\mathbb{N} $, are the time
continuous
interpolations~(\ref{eqtimecont})
of the numerical approximations~(\ref{eqscheme}).
\end{theorem}

Inequality~(\ref{eqmainresult}) shows that the time continuous tamed
Euler approximations~(\ref{eqtimecont}) converge in the strong
$L^p$-sense with the supremum over the time interval $ [0,T] $ inside
the expectation\vspace*{1pt} to the exact solution of the SDE
(\ref{eqSDE}) with the standard convergence order $ \frac{1}{2} $. For
a lower bound of this type of convergence, the reader is referred to
Theorem 3 in M\"{u}ller-Gronbach~\cite{m02} (see also Hofmann,
M\"{u}ller-Gronbach and Ritter~\cite{hmr00a}).\vadjust{\goodbreak}

While the detailed proof of Theorem~\ref{thmmainresult} is postponed to
Section~\ref{secmainresult}, we now outline the central ideas in the
proof of Theorem~\ref{thmmainresult}. The key difficulty in the proof
of Theorem~\ref{thmmain} is to establish that the tamed Euler
approximations~(\ref{eqscheme}) satisfy the a priori moment bounds
%
%
\begin{equation}
\label{eqschemebound2} \sup_{N \in\mathbb{N}} \sup_{ n
\in\{0,1,\ldots,N\} } \mathbb{E} [ \| Y_n^N\|^p ] < \infty
\end{equation}
for all $ p \in[1,\infty) $ (see Lemma~\ref{cmomEuler} in Section
\ref{secmainresult} for the precise statement of this result). After
having verified~(\ref{eqschemebound2}), Theorem 1.1 can, by exploiting
(\ref{eqsecondorder}), at least in the case $p = 2$, be completed
analogously to Theorem 4.4 in Higham, Mao and Stuart~\cite{hms02} in
which strong convergence of the explicit Euler method under the
assumption of the moment bounds~(\ref{eqschemebound2}) has been
established. However, note that, in contrast to the tamed Euler
approximations~(\ref{eqscheme}), the explicit Euler approximations
(\ref{eqEulerscheme}) fail to satisfy such moment bounds for SDEs with
superlinearly growing coefficients (see~(\ref{eqEulerdiv}) here and
Theorem 1 in~\cite{hjk11} for details). It is quite remarkable that
changing the explicit Euler method by a second order term such as in
(\ref{eqsecondorder}) alters the behavior of the numerical method to
such an extent.

Let us now go into details and sketch the central ideas of our proof of
the moment bounds~(\ref{eqschemebound2}). The key idea here for showing
(\ref{eqschemebound2}) is to introduce appropriate stochastic processes
that dominate the tamed Euler approximations~(\ref{eqscheme}) on
appropriate subevents. More formally, let $ D^N_n \dvtx\Omega
\rightarrow[0,\infty) $, $ n \in\{ 0,1,\ldots, N \} $, $ N \in\mathbb{N}
$, be defined by
%
%
\begin{eqnarray}
\label{eqDefDominator}
D_n^N &:=& (\lambda+\|\xi\|)
\exp\Biggl( \lambda+ \sup_{ u \in\{ 0, 1, \ldots, n \} }
\sum_{k=u}^{n-1} \biggl[ \lambda\|\Delta W_k^N
\|^2\nonumber\\[-8pt]\\[-8pt]
&&\hspace*{66pt}{} + \mathbh{1}_{ \{
\norm{Y_k^N} \geq1 \} } \biggl\langle \frac{ Y_k^N }{ \|Y_k^N\| } , \frac{
\sigma(Y_k^N) }{ \|Y_k^N\| } \Delta W_k^N\biggr\rangle\biggr]\Biggr)\nonumber
\end{eqnarray}
for all $ n \in\{0,1,\ldots,N\} $ and all $ N \in\mathbb{N}$ where $
\lambda\in[1,\infty) $ and $ \Delta W^N_n \dvtx\Omega
\rightarrow\mathbb{R}^m $ are defined through $ \lambda:= ( 1 + 2 c + T
+ \| \mu(0) \| + \| \sigma(0) \| )^4 $ and $ \Delta W^N_n := W_{ {
(n+1) T }/{ N } } - W_{ { n T }/{ N } } $ for all $ n \in\{ 0, 1,
\ldots, N-1 \} $ and all $ N \in\mathbb{N} $. We will refer to $ D_n^N$,
$ n \in\{ 0,1, \ldots, N \} $, $ N \in\mathbb{N} $, as \textit{dominating
stochastic processes}. Appropriate subevents are $\Omega_n^N $, $ n
\in\{ 0,1,\ldots, N \}$, $N \in\mathbb{N}$, given by
%
%
\begin{eqnarray}
\label{eqDefOmegaN}
\Omega_n^N &:=& \Bigl\{\omega\in\Omega\dvtx\sup_{ k
\in\{ 0, 1, \ldots, n-1 \} } D_{k}^N(\omega) \leq N^{ { 1 }/({ 2 c })
} ,\nonumber\\[-8pt]\\[-8pt]
&&\hspace*{49pt} {\sup_{ k \in\{ 0, 1, \ldots, n-1 \} }} \| \Delta W_k^N(\omega) \|
\leq1 \Bigr\}\nonumber
\end{eqnarray}
for all $ n \in\{0,1,\ldots,N\} $ and all $ N \in\mathbb{N}$. The main
step of our proof of the moment bounds~(\ref{eqschemebound2}) will be
to establish the pathwise inequality
%
%
\begin{equation}
\label{eqdominated} \mathbh{1}_{ \Omega_n^N } \| Y_n^N \| \leq D_n^N
\end{equation}
for all $ n \in\{ 0, 1, \ldots, N \} $ and all $ N \in\mathbb{N} $
(see Lemma~\ref{ldominator}).
The next step is then to
obtain the moment bounds
%
%
\begin{equation}
\label{eqdominatormean} \limsup_{ N \rightarrow\infty} \mathbb{E}\Bigl[
{\sup_{ n \in\{ 0, 1, \ldots,N \} }} | D_n^N |^p \Bigr] < \infty
\end{equation}
for all $ p \in[1,\infty) $ for the dominating stochastic processes
(see Lemma~\ref{lembdmomodom}). These follow nicely from Doob's
submartingale inequality (see, e.g., Theorem 11.2 in Klenke
\cite{k08b}), from uniform boundedness of $\sigma(Y_k^N)/\|Y_k^N\|$ on
$\{\|Y_k^N\|\geq1\}$ for all $k\in\{0,1,\ldots,N-1\}$, $N\in\N$ and
from the fact that
%
%
\begin{equation}
\E[ \exp( \ld\| \Delta W_k^N \|^2 +\langle v,\Delta W_k^N\rangle ) ] \leq\exp\biggl(
\frac{2 \lambda T m}{N} +\frac{ T \|v\|^2 }{2N} \biggr)
\end{equation}
for all $v\in\R^m$, $k\in\{0,1,\ldots,N-1\}$ and all $ N \in\mathbb{N}
$ with $N\geq4T\ld$ (see Lemmas
\ref{lemexpsqmomnormal}--\ref{lemsupsupsup} for details). Combining
(\ref{eqdominated}) and~(\ref{eqdominatormean}) shows that
%
%
\begin{equation}
\label{eqtamedestimate} \sup_{ N \in\mathbb{N} } \sup_{ n \in\{ 0, 1,
\ldots, N \} } \mathbb{E}[ \mathbh{1}_{ \Omega_n^N } \| Y_n^N \|^p ] <
\infty
\end{equation}
for all $ p \in[1,\infty) $. For proving~(\ref{eqschemebound2}), it
thus remains to verify that
%
%
\begin{equation}
\label{eqtamedestimate2} \sup_{ N \in\mathbb{N} } \sup_{ n \in\{ 0, 1,
\ldots, N \} } \mathbb{E}\bigl[ \mathbh{1}_{ (\Omega_n^N)^c } \| Y_n^N \|^p
\bigr] < \infty
\end{equation}
for all $ p \in[1,\infty) $. This can be achieved by exploiting that
the probability of $ ( \Omega_N^N )^c $ decays rapidly to zero as $ N $
goes to infinity (see Lemma~\ref{lOmegaN} for details) and by using
that the norm $ \| \frac{ {T}/{N} \cdot\mu( Y_n^N ) }{ 1 +
{T}/{N} \cdot\| \mu( Y_n^N ) \| } \| \leq1 $ of the drift term in
(\ref{eqscheme}) is bounded by $1$ due to the taming factor $ 1 / ( 1 +
\frac{T}{N} \cdot\| \mu( Y_n^N ) \| ) $ for $ n \in\{ 0, 1, \ldots, N
\} $ and $ N \in\mathbb{N} $. This is, in fact, the only argument in
our proof of Theorem~\ref{thmmainresult} for which the taming factor in
(\ref{eqscheme}) is needed.

It remains to motivate the pathwise inequality~(\ref{eqdominated}). The
following estimate might at least give an intuition for the case of
globally Lipschitz continuous coefficients. Consider the tamed Euler
scheme for a geometric Brownian motion [SDE~(\ref{eqSDE}) with $d=m=1$,
$\mu(x)=0$ and $\sigma(x)=x$ for all $x\in\R$]. Using $1+x\leq e^x$ for
$x\in\R$ yields that
%
%
\begin{eqnarray}
\|Y_{n+1}^N\|^2 &=& \|Y_{n}^N+Y_{n}^N\Delta W_n^N \|^2
\nonumber\\
&=& \|Y_{n}^N\|^2 +\|Y_{n}^N\Delta W_n^N \|^2 +2\langle Y_n^N,Y_n^N\Delta
W_n^N\rangle
\nonumber\\[-8pt]\\[-8pt]
&=& \|Y_{n}^N\|^2 \biggl( 1+ \|\Delta W_n^N \|^2 +2\biggl\langle\frac{Y_n^N}{\|Y_n^N\|},
\frac{Y_n^N}{\|Y_n^N\|}\Delta W_n^N
\biggr\rangle\biggr)
\nonumber\\
&\leq& \|Y_{n}^N\|^2 \exp\biggl( \|\Delta W_n^N \|^2
+2\biggl\langle\frac{Y_n^N}{\|Y_n^N\|}, \frac{Y_n^N}{\|Y_n^N\|}\Delta W_n^N
\biggr\rangle\biggr)\nonumber
\end{eqnarray}
for all $n\in\{0,1,\ldots,N-1\}$ and all $N\in\N$. Iterating this
inequality leads to an inequality such as~(\ref{eqdominated}). Of
course, the case of superlinearly growing drift coefficients is more
subtle and we refer to Section~\ref{secldominator} for the detailed
proof of inequality~(\ref{eqdominated}).\vadjust{\goodbreak}

Having sketched the central ideas for the proof of our main result, we
now compare Theorem~\ref{thmmainresult} with some related results in
the literature. A priori moment bounds of the form
(\ref{eqschemebound2}) for numerical methods for SDEs with nonglobally
Lipschitz continuous coefficients have been intensively studied in the
literature in recent years.
Whereas the references~\cite{hms02,h96,s10,sm10b,sm10a,hmps10}
deal with a priori moment bounds
of the form~(\ref{eqschemebound2})
for implicit methods,
references
for explicit methods are infrequent. In particular, in 2002, Higham,
Mao and Stuart formulated in~\cite{hms02}, page 1060, the following
open problem: ``In general, it is not clear when such moment bounds
can be expected to hold for explicit methods with $ {f, g \in C^1 }
$'' (drift and diffusion coefficients are denoted by $ f $ and $ g $
in~\cite{hms02} instead of $ \mu$ and $ \sigma$ here). In 2010 an a
priori moment bound for an explicit method has been proved in
Bou-Rabee, Hairer and Vanden-Eijnden~\cite{bhv10}. More precisely, in
the setting of the Langevin equation (see Section 2 in~\cite{bhv10} for
the precise assumptions), Lemma~3.5 in~\cite{bhv10} establishes
exponential moment bounds for a version of the metropolis-adjusted
Langevin algorithm (MALA) with reflection on the boundaries of certain
compact sets. In this paper we concentrate on the explicit numerical
method~(\ref{eqscheme}) in the setting of the SDE~(\ref{eqSDE}). To be
more precise, in Lemma~\ref{cmomEuler} in Section~\ref{secmainresult}
below we show that the explicit method~(\ref{eqscheme}) fulfills
(\ref{eqschemebound2}) for all $ p \in[1,\infty) $.
Furthermore, we observe that one way of deriving the tamed Euler
meth\-od~(\ref{eqscheme}) is to approximate the drift coefficent $ \mu(v) $, $ v
\in\mathbb{R}^d $, in the SDE~(\ref{eqSDE}) by the modified drift
coefficients $ \frac{ \mu(v) }{ 1 + \lambda\| \mu( v ) \| } $, $ v
\in\mathbb{R}^d $, for small $ \lambda\in(0,\infty) $ and then\vspace*{1pt} to apply
the explicit Euler method to these modified SDEs. Approximations of
this type have been used in the literature in order to construct
solutions of nonlinear parabolic unilateral problems (see,
e.g., Palmeri~\cite{Palmeri2000}, Section 3).
Moreover, in the setting of the Langevin equation, Roberts and Tweedie
suggested in 1996 a similar approximation step to~(\ref{eqscheme}) as a
proposal for the Metropolis--Hastings method in order to sample from
the invariant measure of the Langevin SDE (see~\cite{rt96}, Subsection
1.4.3, (12)--(13)). The resulting Metropolis-adjusted method has been
named Metropolis-adjusted Langevin truncated algorithm (MALTA).
Finally, in 1998 a related class of numerical methods has been
considered in Milstein, Platen and Schurz~\cite{mps98} (see
\cite{mps98}, (3.2)--(3.3)) in the case of globally Lipschitz
continuous coefficients of the SDE. To sum up, in the general setting
of the SDE~(\ref{eqSDE}), the tamed Euler method~(\ref{eqscheme}) is,
to the best of our knowledge, the first explicit numerical method that
has been shown to converge strongly to the exact solution of the SDE
(\ref{eqSDE}).

%
%
%
\begin{figure}

\includegraphics{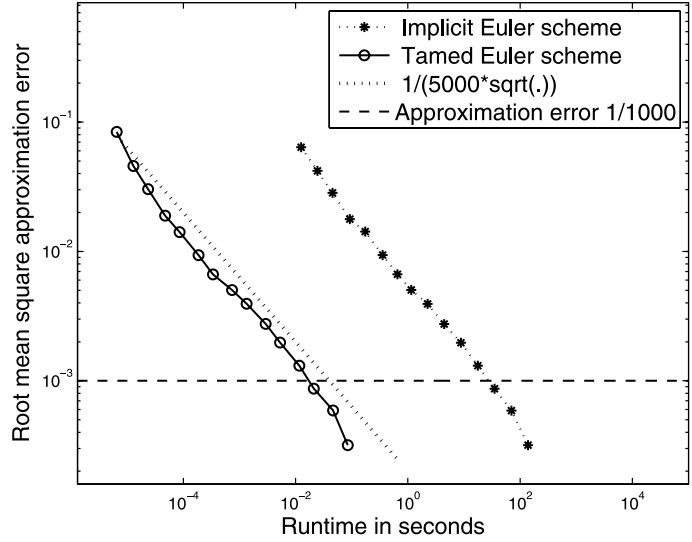}

\caption{Root mean square approximation error $ ( \mathbb{E}\| X_T -
\ttilde{{Y}}_N^N \|^2 )^{ {1}/{2} } $ of the exact solution of
the SDE (\protect\ref{eqGinzburgLandau}) and of the implicit Euler
scheme (\protect\ref{eqimplicitEulerscheme}) and root mean square
approximation error $ ( \mathbb{E}\| X_T - Y_N^N \|^2 )^{ {1}/{2} }
$ of the exact solution of the SDE (\protect\ref{eqGinzburgLandau}) and
of the tamed Euler scheme (\protect\ref{eqscheme}) as function of the
runtime when $ N \in\{ 2^4, 2^5, \ldots, 2^{18} \} $.} \label{fig1}
\end{figure}
%
%

We now compare simulations of the implicit Euler scheme
(\ref{eqimplicitEulerscheme}) and of the tamed Euler scheme
(\ref{eqscheme}). For this we choose $ T = d = m = 1 $, $ \xi= 1 $, $
\mu( x ) = - x^5 $ and $ \sigma(x) = x $ for all $ x \in\mathbb{R} $.
The SDE~(\ref{eqSDE}) thus reads as
%
%
\begin{equation}
\label{eqGinzburgLandau}
dX_t = - X_t^5 \,dt + X_t \,dW_t,\qquad
X_0 = 1,
\end{equation}
for $ t \in[0,1] $. Suppose that the strong approximation problem
(\ref{eqstrongproblem}) of the SDE~(\ref{eqGinzburgLandau}) should be
solved with\vspace*{1pt} the precision of say three decimals, that is, with
precision $ \varepsilon= \frac{1}{1000} $ in~(\ref{eqstrongproblem}).
Figure~\ref{fig1} depicts\vspace*{0.5pt} the root mean square approximation error $ (
\mathbb{E}\| X_T - \ttilde{{Y}}_N^N \|^2 )^{ {1}/{2} } $ of
the exact solution of the SDE~(\ref{eqGinzburgLandau}) and of the
implicit Euler scheme~(\ref{eqimplicitEulerscheme}) and the root mean
square approximation error $ ( \mathbb{E}\| X_T - Y_N^N \|^2 )^{
{1}/{2} } $ of the exact solution of the SDE
(\ref{eqGinzburgLandau}) and of the tamed Euler scheme~(\ref{eqscheme})
as function of the runtime when $N\in\{2^4,2^5,\ldots,2^{18}\}$.
%
%
The zero of the nonlinear equation
that has to be determined in each
time step of the implicit Euler scheme~(\ref{eqimplicitEulerscheme})
is computed approximatively through the function $\mbox{fzero}(\ldots)$
in \textsc{Matlab}.
It turns out that $ \ttilde{{Y}}_{2^{16}}^{2^{16}} $ in the case
of the implicit Euler scheme~(\ref{eqimplicitEulerscheme}) and that $
Y_{2^{16}}^{2^{16}} $ in the case of the tamed Euler
scheme~(\ref{eqscheme}) achieve the precision $ \varepsilon= \frac{1}{1000} $
in~(\ref{eqstrongproblem}).
Following is our \textsc{Matlab} code for
simulating the implicit Euler approximation $
\ttilde{{Y}}^{2^{16}}_{2^{16}} $ [see~(\ref{eqimplicitEulerscheme})]
for the SDE~(\ref{eqGinzburgLandau}):\vspace*{6pt}
\begin{verbatim}
Y = 1; N = 2^16;
for n=1:N
  v = Y + Y*randn/sqrt(N);
  Y = fzero(@(x)x + x^5/N - v, Y);
end
\end{verbatim}\vfill\eject

\noindent Next we specify our \textsc{Matlab} code for calculating the tamed Euler
approximation $ Y^{2^{16}}_{2^{16}} $ [see~(\ref{eqscheme})] for the
SDE~(\ref{eqGinzburgLandau}):\vspace*{6pt}
\begin{verbatim}
Y = 1; N = 2^16;
for n=1:N
  v = -Y^5/N;
  Y = Y + v/(1+abs(v)) + Y*randn/sqrt(N);
end
\end{verbatim}

\vspace*{6pt}
\noindent The above\vspace*{1pt} \textsc{Matlab} code
for calculating the implicit Euler
approximation~$ \ttilde{{Y}}^{2^{16}}_{2^{16}} $
requires, on our computer
running at $ 2.33 $ GHz, about $ 35.1 $ seconds while the above \textsc{Matlab} code for calculating the tamed Euler approximation $
Y^{2^{16}}_{2^{16}} $ requires about $ 0.0212 $ seconds to be evaluated
on the same computer. Thus, on the above computer, the tamed Euler
scheme~(\ref{eqscheme}) for the SDE~(\ref{eqGinzburgLandau}) is more
than one thousand times faster than the implicit Euler scheme
(\ref{eqimplicitEulerscheme}) in achieving a precision of three
decimals in~(\ref{eqstrongproblem}).

%

\section{Further examples}
\label{secexamples}
In this section we present further simulations to illustrate the
efficiency of the tamed Euler scheme~(\ref{eqscheme}). The next example
is a one-dimensional stochastic Ginzburg--Landau equation with
multiplicative noise (see, e.g., Kloeden and Platen
\cite{kp92}, equation (4.52)). More formally,
let $ T = d = m = 1 $,
$ \xi= 1 $,
$ \mu( x ) = x-x^3 $
and $ \sigma( x ) = x $
for all $ x \in\mathbb{R} $.
The SDE~(\ref{eqSDE})
thus reads as
%
%
\begin{equation}
\label{eqexample2}
dX_t = ( X_t - X_t^3 ) \,dt + X_t \,dW_t,\qquad
X_0 = 1,
\end{equation}
for $ t \in[0,1] $.
%
We now use the implicit Euler scheme~(\ref{eqimplicitEulerscheme}) and
the tamed Euler scheme~(\ref{eqscheme}) for approximating the SDE
(\ref{eqexample2}) and we assume again that the strong approximation
problem~(\ref{eqstrongproblem}) of the SDE~(\ref{eqexample2}) should be
solved with the precision of three decimals, that is, with precision $
\varepsilon= \frac{1}{1000} $ in~(\ref{eqstrongproblem}). For
implementing the implicit Euler scheme~(\ref{eqimplicitEulerscheme}) in
the case of the SDE~(\ref{eqexample2}), we observe that the drift
coefficient of the SDE~(\ref{eqexample2}) is a one-dimensional
polynomial of degree three. Roots of one-dimensional polynomials of
degree three are known explicitly thanks to Cardano's method. This
explicit knowledge results in a faster implementation of the implicit
Euler scheme~(\ref{eqimplicitEulerscheme}) than using\vspace*{1pt} the
\textsc{Matlab} function $\mbox{fzero}(\ldots)$.
More precisely, if $p,q\in\mathbb{R}$ and if $ q^2+(p/3)^3 \geq0$, then
the\vspace*{1pt} only real-valued root of $x^3 +px-2q=0$, $x\in\mathbb{R}$, is
$x=(D+q)^{1/3}-(D-q)^{(1/3)}$ where
$D=\sqrt{q^2+(p/3)^3}\in[0,\infty)$. Thus the implicit Euler scheme
(\ref{eqimplicitEulerscheme}) for the SDE~(\ref{eqexample2}) becomes an
explicit scheme and satisfies\looseness=-1
%
%
\begin{eqnarray}
\label{eqimplicitEulerscheme2}
\ttilde{{Y}}_{ n + 1 }^N &=& \bigl( \sqrt{
( q^N_n )^2 + ( N - 1 )^3 /27
} + q^N_n \bigr)^{ {1}/{3} } \nonumber\\[-8pt]\\[-8pt]
&&{}- \bigl( \sqrt{ ( q^N_n )^2 + ( N - 1 )^3 /27
} - q^N_n \bigr)^{ {1}/{3} }\nonumber
\end{eqnarray}\looseness=0
for\vspace*{-1pt} all $ n \in\{ 0, 1, \ldots, N - 1 \} $ and all $ N \in\mathbb{N} $
where $ q^N_n \dvtx\Omega\rightarrow\mathbb{R} $ is defined by $ q^N_n
:= \ttilde{{Y}}_n^N \cdot\frac{N}{2} \cdot( 1 + W_{ { (n +
1) T }/{ N } } - W_{ { n T }/{ N } } ) $ for every\vadjust{\goodbreak} $ n \in\{ 0, 1,
\ldots, N-1 \} $ and every $ N \in\mathbb{N} $.
%
%
\begin{figure}

\includegraphics{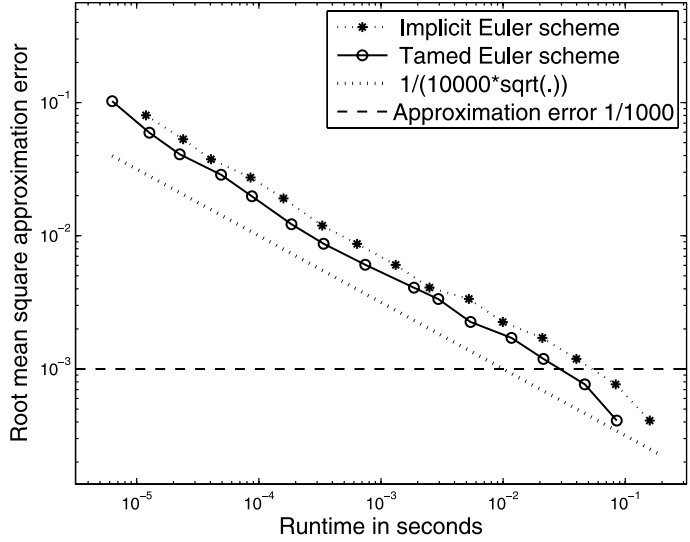}

\caption{Root mean square approximation error $ ( \mathbb{E}\| X_T -
\ttilde{{Y}}_N^N \|^2 )^{ {1}/{2} } $ of the exact solution of
the SDE (\protect\ref{eqexample2}) and of the implicit Euler scheme
(\protect\ref{eqimplicitEulerscheme2}) and root mean square
approximation error $ ( \mathbb{E}\| X_T - Y_N^N \|^2 )^{ {1}/{2} }
$ of the exact solution of the SDE (\protect\ref{eqexample2}) and of
the tamed Euler scheme (\protect\ref{eqscheme}) as function of the
runtime when $N\in\{2^4,2^5,\ldots,2^{18}\}$. } \label{fig2}
\end{figure}
Figure~\ref{fig2} displays the root mean square approximation error $ (
\mathbb{E}\| X_T - \ttilde{{Y}}_N^N \|^2 )^{ {1}/{2} } $ of
the exact solution of the SDE~(\ref{eqexample2}) and of the implicit
Euler scheme~(\ref{eqimplicitEulerscheme2}) and the root mean square
approximation error $ ( \mathbb{E}\| X_T - Y_N^N \|^2 )^{ {1}/{2} }
$ of the exact solution of the SDE~(\ref{eqexample2}) and of the tamed
Euler scheme~(\ref{eqscheme}) as function of the runtime when
$N\in\{2^4,2^5,\ldots,2^{18}\}$. Comparing Figure~\ref{fig1} with
Figure~\ref{fig2} confirms that using explicit knowledge of the roots
of the involved implicit equation results in a much faster
implementation of the implicit Euler scheme. The tamed Euler scheme,
however, is still faster as its implementation does not require the
arithmetical operations for calculating roots. More precisely, it turns
out that $ \ttilde{{Y}}_{2^{17}}^{2^{17}} $ in the case of the
implicit Euler scheme~(\ref{eqimplicitEulerscheme2}) and that $
Y_{2^{17}}^{2^{17}} $ in the case of the tamed Euler scheme
(\ref{eqscheme}) achieve the desired precision $ \varepsilon=
\frac{1}{1000} $ in~(\ref{eqstrongproblem}).
Following is our \textsc{Matlab} code for simulating the implicit Euler
approximation $ \ttilde{{Y}}^{2^{17}}_{2^{17}} $ for the SDE
(\ref{eqexample2}) [see~(\ref{eqimplicitEulerscheme2})]:\vspace*{6pt}
\begin{verbatim}
Y = 1; N = 2^17; v = (N-1)^3/27;
for n=1:N
  q = Y*N*(1+randn/sqrt(N))/2;
  D = sqrt(q^2+v);
  Y = (D+q)^(1/3) - (D-q)^(1/3);
end
\end{verbatim}
\noindent
Next we specify our \textsc{Matlab} code for calculating the tamed Euler
approximation $ Y^{2^{17}}_{2^{17}} $ [see~(\ref{eqscheme})] for the
SDE~(\ref{eqexample2}):\vspace*{6pt}
\begin{verbatim}
Y = 1; N = 2^17;
for n=1:N
  v = (Y-Y^3)/N;
  Y = Y + v/(1+abs(v)) + Y*randn/sqrt(N);
end
\end{verbatim}

\vspace*{6pt} \noindent The above \textsc{Matlab} code for
calculating\vspace*{1pt} the implicit Euler approximation~$
\ttilde{{Y}}^{2^{17}}_{2^{17}} $ requires, on our computer running at $
2.33 $ GHz, about $ 0.0836 $ seconds while the above \textsc{Matlab}
code for calculating the tamed Euler approximation~$
Y^{2^{17}}_{2^{17}} $ requires about $ 0.0467 $ seconds to be evaluated
on the same computer. Thus the tamed Euler scheme~(\ref{eqscheme}) is
on our computer even in the case of the SDE~(\ref{eqexample2}), where
the implicit Euler scheme can be computed explicitly, almost two times
faster than the implicit Euler scheme~(\ref{eqimplicitEulerscheme2}) in
achieving a~precision of three decimals in~(\ref{eqstrongproblem}).

Our last example is a multi-dimensional Langevin equation. More
precisely, we consider the motion of a Brownian particle of unit mass
in the $d$-dimensional potential $ \frac{1}{4} \|x\|^4 - \frac{1}{2} \|
x \|^2 $, $ x \in\mathbb{R}^d $, where $d\in\N$. The corresponding
force on the particle is then $ x - \|x\|^2 \cdot x $, $x\in\R^d$. More
formally, let $T=1$, $m=d \in\mathbb{N}$, $\xi=(0, 0, \ldots, 0 ) $,
$\mu(x)=x-\|x\|^2\cdot x$ and let $\sigma(x)=\mathrm{I}$ be the
identity matrix for all $x\in\mathbb{R}^d$. Thus the SDE~(\ref{eqSDE})
reduces to the Langevin equation
%
%
\begin{equation}
\label{eqexample4}
dX_t = ( X_t - \|X_t\|^2\cdot X_t ) \,dt + dW_t,\qquad
X_0 = 0,
\end{equation}
for $ t \in[0,1] $.
Here
$W\dvtx[0,1]\times\Omega\to\R^d$
is a $d$-dimensional standard Brownian
motion.
%
%
\begin{figure}

\includegraphics{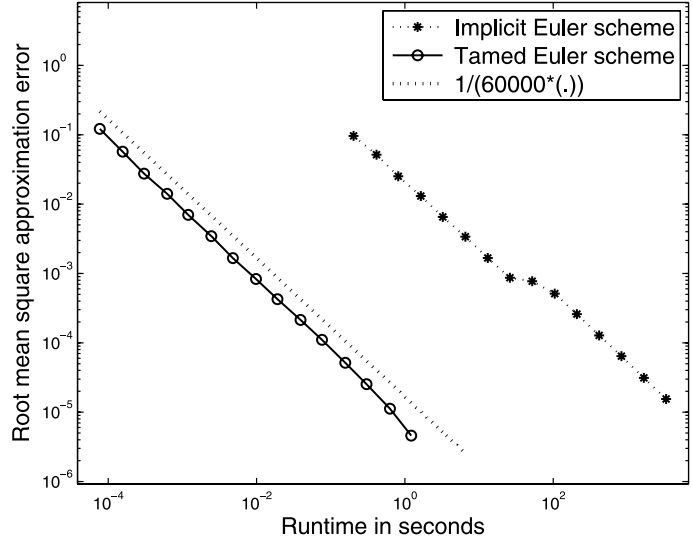}

\caption{Root mean square approximation error $ ( \mathbb{E}\| X_T -
\ttilde{{Y}}_N^N \|^2 )^{ {1}/{2} } $ of the exact solution of
the SDE (\protect\ref{eqexample4}) with $d=10$ and of the implicit
Euler scheme (\protect\ref{eqimplicitEulerscheme}) and root mean square
approximation error $ ( \mathbb{E}\| X_T - Y_N^N \|^2 )^{ {1}/{2} }
$ of the exact solution of the SDE (\protect\ref{eqexample4}) with
$d=10$ and of the tamed Euler scheme (\protect\ref{eqscheme}) as
function of the runtime when $N\in\{2^4,2^5,\ldots,2^{18}\}$.}
\label{fig4}
\end{figure}
%
%
Figure~\ref{fig4} displays the root mean square approximation error $ (
\mathbb{E}\| X_T - \ttilde{{Y}}_N^N \|^2 )^{ {1}/{2} } $ of
the exact solution of the SDE~(\ref{eqexample4}) with $d=10$ and of the
implicit Euler scheme~(\ref{eqimplicitEulerscheme}) and the root mean
square approximation error $ ( \mathbb{E}\| X_T - Y_N^N \|^2 )^{
{1}/{2} } $ of the exact solution of the SDE~(\ref{eqexample4})
with $d=10$ and of the tamed Euler scheme~(\ref{eqscheme}) as function
of the runtime when $N\in\{2^4,2^5,\ldots,2^{18}\}$. We see that both
numerical approximations of the SDE~(\ref{eqexample4}) apparantly
converge with rate $1$. This is presumably due to the additive noise in
(\ref{eqexample4}).
Note that we used the \textsc{Matlab} function $\mbox{fsolve}(\ldots)$ in
our implementation of the implicit Euler scheme
(\ref{eqimplicitEulerscheme}) for the SDE~(\ref{eqexample4}) as the
\textsc{Matlab} function $\mbox{fzero}(\ldots)$ used for the numerical
simulations in Figure~\ref{fig1} is restricted to one dimension
($d=1$). 


Some applications involve high-dimensional SDEs (see, e.g., Beskos and
Stuart~\cite{bs09}, Section 2.1). Then the above implementation of the
implicit Euler method has an additional disadvantage. The \textsc{Matlab}
command $\mbox{fsolve}(\ldots)$ uses (by default) the
``trust-region-dogleg'' algorithm which calculates Jacobian matrices.
Thus the computational effort increases quadratically with the
dimension $ d \in\mathbb{N} $ for every fixed $N\in\mathbb{N}$. To
visualize this we have plotted the runtime of the calculation of the
implicit Euler approximation $\ttilde{{Y}}_{128}^{128}$ of the SDE
(\ref{eqexample4}) as a function of the dimension
$d\in\{10,20,30,\ldots,400\}$.
%
%
\begin{figure}

\includegraphics{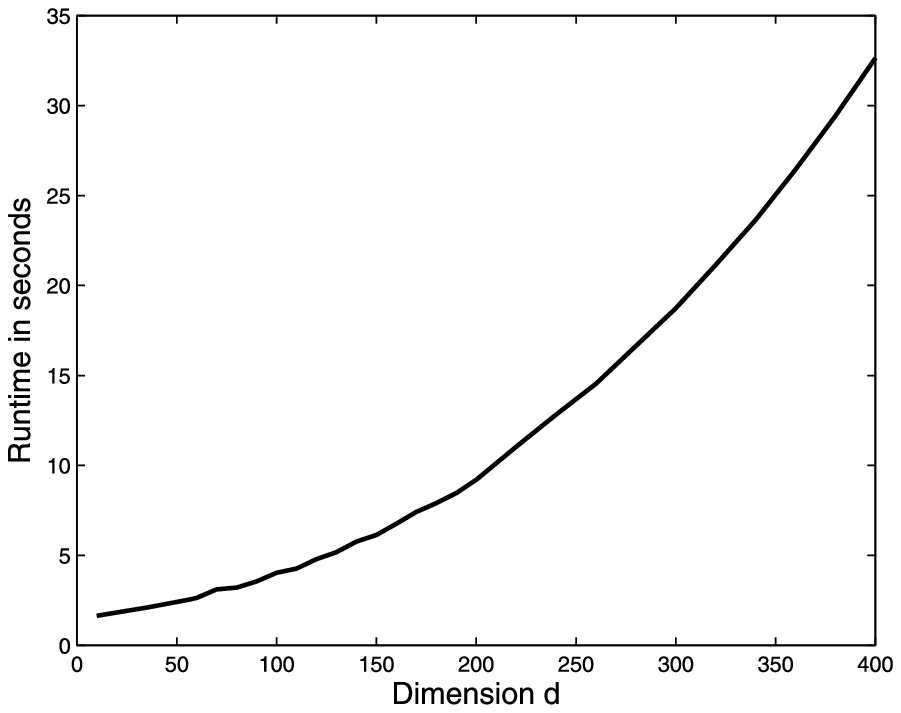}
\caption{Runtime of the calculation of the implicit Euler approximation
$\ttilde{{Y}}_{128}^{128}$ of the SDE~(\protect\ref{eqexample4})
as a function of the dimension $d\in\{10,20,30,\ldots,400\}$.}
\label{fdimImplicit}
\end{figure}
Figure~\ref{fdimImplicit} suggests a quadratic dependence of the
runtime of the implicit Euler method on the dimension in case of the
SDE~(\ref{eqexample4}). In contrast, the tamed Euler method is linear
in the dimension (except that the evaluation of the coefficient
functions might increase quadratically in the dimension).
%
%
\begin{figure}

\includegraphics{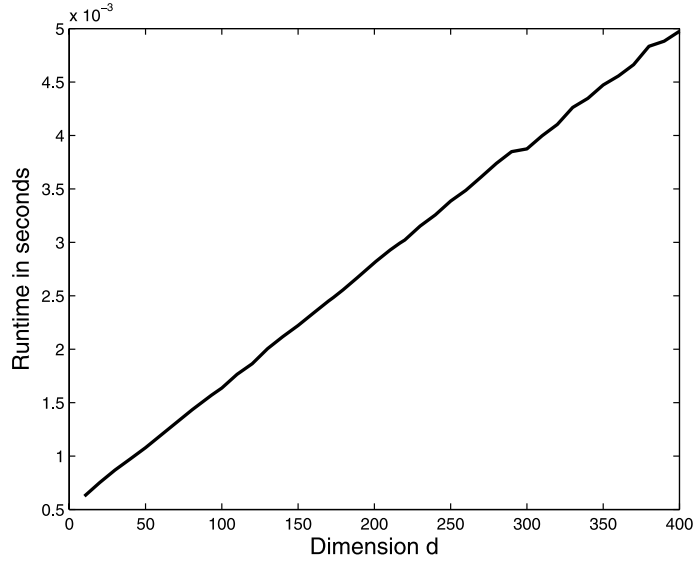}
\caption{Runtime of the calculation of the tamed Euler approximation
$Y_{128}^{128}$ of the SDE~(\protect\ref{eqexample4}) as a function of
the dimension $d\in\{10,20,30,\ldots,400\}$.}
\label{fdimTamed}
\end{figure}
Figure~\ref{fdimTamed} suggests a linear dependence of the runtime of
the tamed Euler method on the dimension $d\in\{10,20,30,\ldots,400\}$
in case of the SDE~(\ref{eqexample4}).

\section{\texorpdfstring{Proof of Theorem \protect\ref{thmmain}}{Proof of Theorem 1.1}}
\label{secmainresult}
In order to simplify
the notation introduced in Section~\ref{secintro},
we introduce
the mappings
$ \alpha^N_n \dvtx\Omega
\rightarrow\mathbb{R} $
defined by
%
%
\begin{equation}
\label{eqDefAlpha}
\alpha_n^N := \mathbh{1}_{ \{ \| Y_n^N \| \geq1 \} }
\biggl\langle \frac{ Y_n^N }{ \| Y_n^N \| }, \frac{ \sigma( Y_n^N ) }{ \| Y_n^N \|
} \Delta W_n^N
\biggr\rangle
\end{equation}
for all $ n \in\{ 0, 1, \ldots, N-1 \} $ and all $ N \in\mathbb{N} $.
Using this notation, the dominating stochastic processes $ D^N_n $, $ n
\in\{ 0,1, \ldots, N \} $, $ N \in\mathbb{N} $, [see
(\ref{eqDefDominator})] simplify to
\[
D_n^N = ( \lambda+ \|\xi\| ) \cdot\exp\Biggl( \lambda+ \sup_{ u \in\{ 0, 1,
\ldots, n \} } \sum_{k=u}^{n-1} [ \lambda\|\Delta W_k^N \|^2 +
\alpha^N_k ] \Biggr)
\]
for all $ n \in\{ 0, 1, \ldots, N \} $ and all $ N \in\mathbb{N} $.
Moreover, we denote by $ \vec{e}_1 = (1, 0, \ldots,\break 0) \in\mathbb{R}^m
$, $ \vec{e}_2 = (0, 1, 0, \ldots, 0) \in\mathbb{R}^m , \ldots,
\vec{e}_m = (0, \ldots, 0, 1) \in\mathbb{R}^m $ the unit vectors in the
$ \mathbb{R}^m $. Additionally, we use the mappings $ \sigma_i
\dvtx\mathbb{R}^d \rightarrow\mathbb{R}^d $, $ i \in\{ 1,2, \ldots,\break m \}
$, given by $ \sigma_i(x) = ( \sigma_{1,i}(x), \ldots, \sigma_{d,i}(x) )
= \sigma(x) \vec{e}_i $ for all $ x \in\mathbb{R}^d $ and all $ i \in\{
1,2, \ldots, m \} $. The SDE~(\ref{eqSDEsol}) can thus be written as
%
%
\begin{equation}
X_t = \xi+ \int_0^t \mu( X_s ) \,ds + \sum_{ i=1 }^m \int_0^t
\sigma_i( X_s ) \,dW^{(i)}_s
\end{equation}
for all $ t \in[0,T] $ $ \mathbb{P} $-a.s. Our proof of Theorem
\ref{thmmainresult} relies on the following lemmas.
%
%
\begin{lemma}[(Dominator lemma)] \label{ldominator}
Let $Y_n^N\dvtx\Omega\to\mathbb{R}^d$, $ D^N_n\dvtx\Omega\to[0,\infty)
$ and $ \Omega^N_n \in\mathcal{F} $ for $n \in\{0,1,\ldots,N\}$ and $ N
\in\mathbb{N} $ be given by~(\ref{eqscheme}),~(\ref{eqDefDominator})
and~(\ref{eqDefOmegaN}). Then we have that
%
%
\begin{equation}
\label{eqdominator}
\mathbh{1}_{ \Omega_n^N } \|Y_n^N \| \leq D_n^N
\end{equation}
for all $n\in\{0,1,\ldots,N\}$
and all $N\in\N$.
\end{lemma}
%
%
\begin{lemma}
\label{lemexpsqmomnormal} Let $ n \in\mathbb{N} $ and let $Z
\dvtx\Omega\rightarrow\mathbb{R}^n $ be an $n$-dimensional standard
normal random variable. Then we have that
%
%
\begin{equation}
\label{eqexpsqmomnormal} \mathbb{E}[ \exp( a \| Z \|^2 ) ] = ( 1 - 2a
)^{ -{n}/{2} } \leq e^{ 2 a n }
\end{equation}
for all $ a \in[ 0, \frac{1}{4} ]$.
\end{lemma}
%
%
\begin{lemma}
\label{lemexplemma}
We have that
%
%
\begin{equation}
\mathop{\sup_{N \in\mathbb{N},}}_{N \geq4 \lambda p T} \mathbb{E}\Biggl[ \exp
\Biggl( p \lambda\sum_{ k = 0 }^{ N - 1 } \| \Delta W_k^N \|^2 \Biggr) \Biggr] < \infty
\end{equation}
for all $ p \in[1,\infty) $.
\end{lemma}
%
%
%
%
\begin{lemma}
\label{lemsupsupsup} Let $ \alpha^N_n \dvtx\Omega\to\mathbb{R} $ for $
n \in\{ 0, 1, \ldots, N \} $ and $ N \in\mathbb{N} $ be given by
(\ref{eqDefAlpha}). Then we have that
%
%
\begin{equation}
\sup_{ z \in\{ -1, 1 \} } \sup_{ N \in\mathbb{N} } \mathbb{E} \Biggl[ \sup_{
n \in\{ 0, 1, \ldots, N \} } \exp\Biggl( p z \sum_{ k = 0 }^{ n - 1 }
\alpha_k^N \Biggr) \Biggr] < \infty
\end{equation}
for all $ p \in[1, \infty) $.
\end{lemma}
%
%
\begin{lemma}[(Uniformly
bounded moments of the dominating stochastic processes)]
\label{lembdmomodom} Let $D_n^N\dvtx\Omega\to[0,\infty)$ for $n
\in\{0,1,\ldots,N\}$ and $ N \in\mathbb{N} $ be given by
(\ref{eqDefDominator}). Then we have that
%
%
\begin{equation}
\mathop{\sup_{N \in\mathbb{N},}}_{N \geq8 \lambda p T} \mathbb{E}\Bigl[
{\sup_{ n \in\{ 0, 1, \ldots, N \} }} | D_n^N |^p \Bigr] < \infty
\end{equation}
for all $ p \in[1, \infty) $.
\end{lemma}
%
%
%
%
\begin{lemma}[(Estimation of
the probability of the complement of $ \Omega^N_N $ for \mbox{$ N
\in\mathbb{N} $})] \label{lOmegaN} Let $ \Omega_N^N \in\mathcal{F} $ for
$ N \in\mathbb{N} $ be given by~(\ref{eqDefOmegaN}). Then we have that
%
%
\begin{equation}
\label{eqOmegaNgoal} \sup_{ N \in\N} \bigl( N^p \cdot\mathbb{P}[ (
\Omega_N^N )^c ] \bigr) < \infty
\end{equation}
for all $p\in[1,\infty)$.
\end{lemma}
%
%
\begin{lemma}[(Time continuous Burkholder--Davis--Gundy type inequality)]
\label{lBDG} Let $k\in\N$ and let $Z
\dvtx[0,T]\times\Omega\to\R^{k \times m}$ be a predictable stochastic
process satisfying $ \mathbb{P}[ \int_0^T \| Z_s \|^2 \,ds < \infty] = 1
$. Then we obtain that
%
%
\begin{equation}\qquad
\biggl\|\sup_{s\in[0,t]} \biggl\| \int_0^s Z_u \,dW_u \biggr\|
\biggr\|_{L^p(\Omega;\R)} \leq p \Biggl( \int_0^t \sum_{i=1}^m \| Z_s
\vec{e}_i \|_{L^p(\Omega,\R^{k})}^2 \,ds \Biggr)^{ {1}/{2}}
\end{equation}
for all $t\in[0,T]$ and
all $p\in[2,\infty)$.
\end{lemma}
%
%
\begin{lemma}[(Time discrete
Burkholder--Davis--Gundy type inequality)] \label{cBDG} Let
$k\in\N$ and let $Z_l^N \dvtx\Omega\to\R^{k \times m}$, $ l \in\{ 0, 1,
\ldots, N-1 \} $, $ N \in\mathbb{N} $, be a family of mappings such that
$Z_l^N \dvtx\Omega\to\R^{k \times m}$ is $ \mathcal{F}_{{lT}/{N}}
/ \mathcal{B}( \mathbb{R}^{k \times m} ) $-measurable for all
$l\in\{0,1,\ldots,N-1\}$ and all $N\in\N$. Then we obtain that
%
%
\begin{equation}
\Biggl\| \sup_{j\in\{0,1,\ldots,n\}} \Biggl\| \sum_{l=0}^{j-1}Z_l^N\Delta W_l^N \Biggr\|
\Biggr\|_{L^p(\Omega;\R)} \leq p\Biggl(\sum_{l=0}^{n-1}\sum_{i=1}^m \|Z_l^N
\vec{e}_i \|_{L^p(\Omega;\R^{k})}^2 \frac{T}{N} \Biggr)^{{1}/{2}}\hspace*{-35pt}
\end{equation}
for all $n\in\{0,1,\ldots,N\}$, $N\in\N$ and all $p\in[2,\infty)$.
\end{lemma}
%
%
\begin{lemma}[(Uniformly bounded
moments of the tamed Euler approximations)] \label{cmomEuler} Let
$Y_n^N\dvtx\Omega\to\mathbb{R}^d$ for $n \in\{0,1,\ldots,N\}$ and $ N
\in\mathbb{N} $ be given by~(\ref{eqscheme}). Then we have that
%
%
\begin{equation}\label{eqschemebound}
\sup_{N \in\mathbb{N}} \sup_{ n \in\{0,1,\ldots,N\} } \mathbb{E}[ \|
Y_n^N \|^p ] < \infty
\end{equation}
for all $p \in[1,\infty) $.
\end{lemma}
%
%
\begin{lemma}\label{lemmubound}
Let $Y_n^N\dvtx\Omega\to\mathbb{R}^d$ for $n \in\{0,1,\ldots,N\}$ and $
N \in\mathbb{N} $ be given by~(\ref{eqscheme}). Then we have that
\[
\sup_{N \in\mathbb{N}} \sup_{ n \in\{0,1,\ldots,N\} } \mathbb{E}[ \|
\mu( Y_n^N ) \|^p ] < \infty,\qquad
\sup_{N \in\mathbb{N}} \sup_{ n \in\{0,1,\ldots,N\} } \mathbb{E}[ \|
\sigma( Y_n^N ) \|^p ] < \infty
\]
for all $p \in[1,\infty) $.
\end{lemma}

The proofs of
Lemmas~\ref{ldominator}--\ref{lemmubound}
can be found in
Sections~\ref{secldominator}--\ref{secproofmubound}.
%
Using Lemmas~\ref{cmomEuler} and~\ref{lemmubound},
the proof of Theorem~\ref{thmmain} is then
completed in Section~\ref{secproofthm}.
%
%
\subsection{\texorpdfstring{Proof of Lemma \protect\ref{ldominator}}{Proof of Lemma 3.1}}
\label{secldominator}

First of all, note that $\|\Delta W_n^N\|\leq1$ on $\Omega_{n+1}^N$ for
all $n\in\{0,1,\ldots,N-1\}$ and all $ N \in\mathbb{N} $. The global
Lipschitz continuity of $ \sigma$ and the polynomial growth bound on $
\mu' $ (see Section~\ref{secintro}) therefore imply that
%
%
\begin{eqnarray}
\label{eqestimateinunitball}
\| Y_{ n + 1 }^N \| &\leq& \|
Y_n^N \| + \frac{ T }{ N } \| \mu( Y_n^N ) \| + \| \sigma( Y_n^N ) \|
\| \Delta W_n^N \|
\nonumber\\[-1pt]
&\leq&
1 + T \| \mu( Y_n^N ) - \mu( 0 ) \| + T \| \mu(0)
\|\nonumber\\[-8.5pt]\\[-8.5pt]
&&{} + \| \sigma( Y_n^N
) - \sigma(0) \| + \| \sigma(0) \|
\nonumber\\[-1pt]
&\leq&
1 + 2 T c + T \| \mu(0) \| + c + \| \sigma(0) \| \leq\lambda\nonumber
\end{eqnarray}
on $ \Omega_{n+1}^N \cap\{ \omega\in\Omega\dvtx\| Y^N_n(\omega) \|
\leq1 \} $ for all $ n \in\{ 0,1,\ldots, N-1 \}$ and all $ N
\in\mathbb{N} $.

Moreover, the Cauchy--Schwarz inequality and the estimate $ a \cdot b
\leq\frac{ a^2 }{ 2 } + \frac{ b^2 }{ 2 } $ for all $ a, b
\in\mathbb{R} $ show that
%
%
\begin{eqnarray}
\label{eqfirstestimate}
\|Y_{n+1}^N\|^2
&=& \biggl\| Y_n^N + \frac{ ({T}/{N}) \mu(Y_n^N ) }{ 1 + ({T}/{N}) \norm{
\mu(Y_n^N) } } + \sigma(Y_n^N ) \Delta W_n^N \biggr\|^2 \nonumber\\[-1pt]
&=&
\|Y_{n}^N\|^2 +\frac{ \|({T}/{N})\mu(Y_{n}^N) \|^2 }{ (1+({T}/{N})
\norm{\mu(Y_n^N)} )^2 } + \|\sigma(Y_{n}^N) \Delta W_n^N \|^2  \nonumber\\[-1pt]
&&{}+ \frac{
2 \langle Y_n^N ,({T}/{N})\mu(Y_{n}^N ) \rangle}{ 1+({T}/{N})
\norm{ \mu(Y_n^N) } } + 2 \langle Y_n^N ,\sigma(Y_{n}^N )
\Delta W_n^N \rangle\nonumber\\[-1pt]
&&{}+ \frac{ 2 \langle({T}/{N}) \mu(Y_n^N) ,\sigma
(Y_{n}^N) \Delta W_n^N \rangle}{ 1 + ({T}/{N}) \norm{\mu(Y_n^N) } }
\nonumber\\[-1pt]
&\leq& \|Y_{n}^N \|^2 + \frac{T^2}{N^2} \|\mu(Y_{n}^N ) \|^2
+ \|\sigma(Y_{n}^N) \|^2 \|\Delta W_n^N \|^2  \\[-1pt]
&&{} + \frac{ ({2T}/{N}) \langle
Y_n^N, \mu( Y_n^N )
\rangle
}{ 1 + ({T}/{N}) \| \mu( Y_n^N ) \| }+ 2 \langle Y_n^N
,\sigma(Y_{n}^N) \Delta W_n^N \rangle\nonumber\\[-1pt]
&&{}+ \frac{2 T}{N} \|\mu(Y_n^N) \|
\|\sigma(Y_{n}^N) \| \|\Delta W_n^N \| \nonumber\\[-1pt]
&\leq& \|Y_{n}^N \|^2
+ \frac{ 2 T^2 }{ N^2 } \|\mu(Y_{n}^N) \|^2 + 2 \|\sigma(Y_{n}^N) \|^2
\|\Delta W_n^N \|^2  \nonumber\\[-1pt]
&&{} + \frac{ ({2T}/{N}) \langle Y_n^N, \mu( Y_n^N ) \rangle }{ 1
+ ({T}/{N}) \| \mu( Y_n^N ) \| }+2 \langle Y_n^N
,\sigma(Y_{n}^N) \Delta W_n^N \rangle\nonumber
\end{eqnarray}
on $\Omega$ for all $ n \in\{ 0,1, \ldots, N-1 \} $ and all $ N
\in\mathbb{N} $. Additionally, the global Lipschitz continuity of $
\sigma$ (see Section~\ref{secintro}) implies that
%
%
\begin{eqnarray}
\| \sigma(x) \|^2 &\leq& \bigl( \| \sigma(x) - \sigma(0) \| + \| \sigma(0) \|
\bigr)^2 \nonumber\\
&\leq&\bigl( c \| x \| + \| \sigma(0) \| \bigr)^2
\leq \bigl( c + \| \sigma(0) \| \bigr)^2 \| x \|^2 \\
&\leq&\lambda\| x \|^2
\nonumber\vadjust{\goodbreak}
\end{eqnarray}
for all $x\in\R^d$ with $ \norm{x} \geq1 $ and the global one-sided
Lipschitz continuity of~$\mu$ (see Section~\ref{secintro}) gives that
%
%
\begin{eqnarray}
\langle x, \mu(x) \rangle &=& \langle x, \mu(x) -
\mu(0) \rangle + \langle x, \mu(0) \rangle \nonumber\\[0.2pt]
&\leq& c \| x \|^2 + \| x \| \| \mu(0) \|
\leq \bigl( c + \| \mu(0) \| \bigr) \| x \|^2 \\[0.2pt]
&\leq&\sqrt{ \lambda} \| x \|^2
\nonumber
\end{eqnarray}
for all $x\in\R^d$ with $ \norm{x} \geq1 $. Furthermore, the polynomial
growth bound on~$\mu'$ (see Section~\ref{secintro}) yields that
%
%
\begin{eqnarray}
\label{eqlastestimate}
\| \mu(x) \|^2 &\leq&\bigl( \| \mu(x) - \mu(0) \| + \|
\mu(0) \| \bigr)^2 \nonumber\\
&\leq&\bigl( c ( 1 + \| x \|^c ) \| x \| + \| \mu(0) \| \bigr)^2
\leq\bigl( 2 c \| x \|^{ ( c + 1 ) } + \| \mu(0) \| \bigr)^2\\
&\leq& \bigl( 2 c + \|
\mu(0) \| \bigr)^2 \| x \|^{ 2 ( c + 1 ) } \leq N \sqrt{ \lambda} \| x \|^2
\nonumber
\end{eqnarray}
for all $x\in\R^d$ with $ 1 \leq\norm{x} \leq N^{ { 1 }/({ 2c }) } $
and all $ N \in\mathbb{N} $. Combining
(\ref{eqfirstestimate})--(\ref{eqlastestimate}) and $ T^2 + T
\leq\sqrt{ \lambda} $ then gives that
%
%
\begin{eqnarray}\label{eqestimateoutofunitball}
\|Y_{n+1}^N\|^2 &\leq& \|Y_{n}^N\|^2 + \frac{ 2 T^2 \sqrt{
\lambda} }{ N } \|Y_{n}^N\|^2 + 2 \lambda\|Y_{n}^N\|^2 \|\Delta
W_n^N\|^2\nonumber\\
&&{} + \frac{ 2 T \sqrt{ \lambda} }{ N } \|Y^N_n \|^2
+ 2 \langle Y_n^N ,\sigma(Y_{n}^N ) \Delta W_n^N \rangle
\nonumber\\[-8pt]\\[-8pt]
&=&
\|Y_{n}^N \|^2 \biggl( 1 + \frac{ 2 ( T^2 + T
) \sqrt{ \lambda} }{N} + 2 \lambda\|\Delta W_n^N\|^2 + 2 \alpha^N_n \biggr)
\nonumber\\
&\leq&
\|Y_{n}^N \|^2 \exp\biggl( \frac{ 2 \lambda}{ N } + 2
\lambda\|\Delta W_n^N\|^2 + 2 \alpha^N_n \biggr)\nonumber
\end{eqnarray}
on $ \{ \omega\in\Omega\dvtx1 \leq\| Y^N_n( \omega) \| \leq N^{
{ 1 }/({ 2 c }) } \} $ for all $ n \in\{ 0,1,\ldots, N-1 \}$ and all
\mbox{$N \in\mathbb{N} $}.

Additionally, we use the mappings $ \tau^N_l \dvtx\Omega\rightarrow\{
-1, 0, 1, \ldots, l-1 \} $, $ l \in\{ 0,1,\break\ldots, N \}$, $ N
\in\mathbb{N}$, given by
%
%
\begin{equation}
\label{eqtaudef} \tau_{l}^N(\omega) := \max\bigl( \{ -1 \} \cup\bigl\{ n \in\{
0,1,\ldots, l-1 \} | \|Y_n^N(\omega) \| \leq1 \bigr\} \bigr)
\end{equation}
for all $ \omega\in\Omega$, $ l \in\{ 0,1,\ldots, N \}$ and all $ N
\in\mathbb{N}$.

With the estimates~(\ref{eqestimateinunitball}) and
(\ref{eqestimateoutofunitball}) at hand, we now establish
(\ref{eqdominator}) by induction on $n\in\{0,1,\ldots,N\}$ where $ N
\in\mathbb{N} $ is fixed. The base case $n=0$ is trivial. Now let $ l
\in\{0,1,\ldots, N-1\}$ be arbitrary and assume that
inequality~(\ref{eqdominator}) holds for all $n \in\{0,1,\ldots,l\}$. We then show
inequality~(\ref{eqdominator}) for $n=l+1$. More formally, we now
establish that
%
%
\begin{equation} \label{eqdominatoromega}
\|Y_{ l + 1 }^N(\omega) \|
\leq D_{ l + 1 }^N(\omega)
\end{equation}
for all $\omega\in\Omega^N_{l+1}$. To this end let $ \omega
\in\Omega_{l+1}^N$ be arbitrary. From the induction hypothesis and from
$ \omega\in\Omega_{l+1}^N \subseteq\Omega_{n+1}^N $ it follows that $
\|Y_n^N(\omega)\| \leq D_n^N(\omega) \leq N^{ { 1 }/({ 2c }) } $ for
all $n \in\{ 0,1,\ldots,l \}$. By definition~(\ref{eqtaudef}) we
therefore obtain that $ 1 \leq\|Y_n^N(\omega)\| \leq N^{ { 1 }/({ 2c
}) } $ for all $ n \in\{ \tau^N_{l+1}(\omega)+1,
\tau^N_{l+1}(\omega)+2, \ldots, l \} $. Estimate~(\ref{eqestimateoutofunitball}) thus gives that
%
%
\begin{equation}\label{eqiterate}
\| Y^N_{ n+1 }(\omega) \| \leq\|Y_{n}^N(\omega) \|\cdot\exp\biggl( \frac{
\lambda}{ N } + \lambda\|\Delta W_n^N \|^2 + \alpha_n^N(\omega) \biggr)
\end{equation}
for all $ n \in\{ \tau^N_{l+1}(\omega)+1, \tau^N_{l+1}(\omega)+2,
\ldots, l \} $. Iterating~(\ref{eqiterate}) hence yields that
\begin{eqnarray*}
\hspace*{-3pt}&&\| Y_{ l + 1 }^N( \omega) \|
\\
\hspace*{-3pt}&&\quad\leq
\| Y_l^N( \omega) \| \cdot\exp\biggl( \frac{ \lambda}{ N } +
\lambda\| \Delta W_l^N( \omega) \|^2 + \alpha_l^N( \omega) \biggr)
\\
\hspace*{-3pt}&&\quad\leq
\cdots\leq\bigl\| Y_{ \tau_{ l + 1 }^N( \omega) + 1 }^N( \omega)
\bigr\| \cdot\exp\Biggl( \sum_{ n = \tau_{ l + 1 }^N( \omega) + 1 }^l \biggl( \frac{
\lambda}{ N } + \lambda\| \Delta W_n^N( \omega) \|^2 + \alpha_n^N(
\omega) \biggr) \Biggr)
\\
\hspace*{-3pt}&&\quad\leq
\bigl\| Y_{ \tau_{ l + 1 }^N( \omega) + 1 }^N( \omega) \bigr\| \cdot
\exp\Biggl( \lambda+ \sup_{ u \in\{ 0,1,\ldots,l+1 \} } \sum_{ n = u }^l \bigl(
\lambda\| \Delta W_n^N( \omega) \|^2 + \alpha_n^N( \omega) \bigr) \Biggr) .
\end{eqnarray*}
Estimate~(\ref{eqestimateinunitball})
therefore shows that
%
%
\begin{eqnarray}
&&\| Y_{ l + 1 }^N( \omega) \|\hspace*{-25pt}
\nonumber\\
&&\qquad\leq
\bigl( \lambda+ \| \xi(\omega) \| \bigr) \cdot\exp\Biggl( \lambda+ \sup_{
u \in\{ 0,1,\ldots,l+1 \} } \sum_{ n = u }^l \bigl( \lambda\| \Delta W_n^N(
\omega) \|^2 + \alpha_n^N( \omega) \bigr) \Biggr)\hspace*{-25pt}
\\
&&\qquad=
D^N_{l+1}(\omega) .\nonumber\hspace*{-25pt}
\end{eqnarray}
This finishes the induction step and the proof.

%
%
\subsection{\texorpdfstring{Proof of Lemma \protect\ref{lemexpsqmomnormal}}{Proof of Lemma 3.2}}
We establish~(\ref{eqexpsqmomnormal}) in the case $ n = 1 $. The
general case then follows from independence. For the case $ n = 1 $
note that
%
%
\begin{equation}
\label{eqineqfracexp} \frac{1}{ ( 1-x ) } \leq e^{2x}
\end{equation}
for all $ x \in[0,\frac{1}{2}] $.
Inequality~(\ref{eqineqfracexp}) then shows that
\begin{eqnarray*}
\mathbb{E}[ \exp( a Z^2 ) ] &=& \int_{ \mathbb{R} } e^{ a x^2 } \frac{ 1
}{ \sqrt{ 2 \pi} } e^{ - { x^2 }/{ 2 } } \,dx = \int_{ \mathbb{R} }
\frac{ 1 }{ \sqrt{ 2 \pi} } e^{ { -x^2 }/({ 2( 1 - 2a )^{-1} }) }
\,dx\\
&=& \frac{ 1 }{ \sqrt{ 1 - 2a } } \leq e^{ 2a }
\end{eqnarray*}
for all $ a \in[0,\frac{1}{4}] $. This completes the proof of
Lemma~\ref{lemexpsqmomnormal}.

%
\subsection{\texorpdfstring{Proof of Lemma \protect\ref{lemexplemma}}{Proof of Lemma 3.3}}
The independence of the random variables $ \Delta W_k^N $, $ k \in\{ 0,
1, \ldots, N-1 \} $, and Lemma~\ref{lemexpsqmomnormal} show that
%
%
\begin{eqnarray}
\mathbb{E}\Biggl[ \exp\Biggl( \lambda p \sum_{ k = 0 }^{ N - 1 } \| \Delta W_k^N
\|^2 \Biggr) \Biggr]
&=& \prod_{ k = 0 }^{ N - 1 } \mathbb{E}[ \exp( \lambda p \|
\Delta W_k^N \|^2 ) ]\nonumber\\
&=& ( \mathbb{E}[ \exp( \lambda p \|
W_{ { T }/{ N } } \|^2 ) ] )^{ N} \\
&\leq&\exp( 2 \lambda p T m ) <
\infty\nonumber
\end{eqnarray}
for all $ N \in\mathbb{N} \cap[ 4 \lambda p T, \infty) $ and all $ p
\in[ 1, \infty) $.
This completes the proof of
Lem\-ma~\ref{lemexplemma}.

%
\subsection{\texorpdfstring{Proof of Lemma \protect\ref{lemsupsupsup}}{Proof of Lemma 3.4}}

First of all,\vspace*{1pt} note that the time discrete stochastic process $ z \sum_{
k = 0 }^{ n - 1 } \alpha_k^N $, $ n \in\{ 0, 1, \ldots, N \} $, is an
$( \mathcal{F}_{ { n T }/{ N } } )_{ n \in\{ 0, 1, \ldots, N \} }
$-martingale for every $ z \in\{ -1, 1 \} $ and every $ N \in\mathbb{N}
$. In particular, we therefore obtain that the time discrete stochastic
process $ \exp( z \sum_{ k = 0 }^{ n - 1 } \alpha_k^N ) $, $ n \in\{ 0,
1, \ldots, N \} $, is a~positive $( \mathcal{F}_{ { n T }/{ N } }
)_{ n \in\{ 0, 1, \ldots, N \} } $-submartingale for every $ z \in\{
-1, 1 \} $ and every $ N \in\mathbb{N} $.
Doob's maximal
inequality (see, e.g., Klenke~\cite{k08b}, Theorem 11.2)
hence shows that
%
%
\begin{equation}
\label{eqsupproofA} \Biggl\| \sup_{ n \in\{ 0, 1, \ldots, N \} } \exp\Biggl( z
\sum_{ k = 0 }^{ n - 1 } \alpha_k^N \Biggr) \Biggr\|_{ L^p( \Omega; \mathbb{R} ) }
\leq\frac{ p }{ ( p - 1 ) } \Biggl\| \exp\Biggl( z \sum_{ k = 0 }^{ N - 1 }
\alpha_k^N \Biggr) \Biggr\|_{ L^p( \Omega; \mathbb{R} ) }\hspace*{-35pt}
\end{equation}
for all $ N \in\mathbb{N} $, $ p \in( 1, \infty) $ and all $ z
\in\{-1,1\} $. Moreover, we have that
\begin{eqnarray*}
&&\mathbb{E}\biggl[ \biggl| p z \mathbh{1}_{ \{ \| x \| \geq1 \} }
\biggl\langle \frac{ x }{ \| x
\| }, \frac{ \sigma( x ) }{ \| x \| } \Delta W_k^N
\biggr\rangle
\biggr|^2 \biggr] \\
&&\qquad= \frac{ p^2 T }{ N } \mathbh{1}_{ \{ \| x \| \geq1 \} } \frac{
\| ( \sigma( x ) )^{*} x \|^2 }{ \| x \|^4 }
\leq
\frac{ p^2 T }{ N } \mathbh{1}_{ \{ \| x \| \geq1 \} } \frac{ \| (
\sigma( x ) )^{*} \|^2 }{ \| x \|^2 } \\
&&\qquad= \frac{ p^2 T }{ N }
\mathbh{1}_{ \{ \| x \| \geq1 \} } \frac{ \| \sigma( x ) \|^2 }{ \| x
\|^2 } \leq\frac{ p^2 T ( c + \| \sigma(0) \| )^2 }{ N }
\end{eqnarray*}
for all $ x \in\mathbb{R}^d $, $ k \in\{ 0, 1, \ldots, N - 1 \} $, $ N
\in\mathbb{N} $, $ p \in[1, \infty) $ and all $ z \in\{ -1, 1 \} $.
Lemma 5.7 in~\cite{hj09b} therefore gives that
%
%
\begin{equation}
\label{eqsupproofB} \mathbb{E}\biggl[ \exp\biggl( p z \mathbh{1}_{ \{ \| x \|
\geq1 \} } \biggl\langle \frac{ x }{ \| x \| }, \frac{ \sigma( x ) }{ \| x \| }
\Delta W_k^N
\biggr\rangle
\biggr) \biggr] \leq\exp\biggl( \frac{ p^2 T ( c + \| \sigma(0) \| )^2 }{ N } \biggr)\hspace*{-35pt}
\end{equation}
for all $ x \in\mathbb{R}^d $, $ k \in\{ 0, 1, \ldots, N - 1 \} $, $ N
\in\mathbb{N} $, $ p \in[1, \infty) $ and all $ z \in\{ -1, 1 \} $.
Estimate~(\ref{eqsupproofB}), in particular, shows that
\[
\mathbb{E}[ \exp( p z \alpha_k^N ) | \mathcal{F}_{ { k T }/{ N } }
] \leq\exp\biggl( \frac{ p^2 T ( c + \| \sigma(0) \| )^2 }{ N } \biggr)
\]
for all $ k \in\{ 0, 1, \ldots, N - 1 \} $, $ N \in\mathbb{N} $, $ p
\in[1, \infty) $ and all $ z \in\{ -1, 1 \} $. Hence, we obtain that
%
%
\begin{eqnarray}\label{eqsupproofC}
\mathbb{E}\Biggl[ \exp\Biggl( p z \sum_{ k = 0 }^{ N - 1 } \alpha_k^N \Biggr) \Biggr]
&=&
\mathbb{E}\Biggl[ \exp\Biggl( p z \sum_{ k = 0 }^{ N - 2 } \alpha_k^N \Biggr)
\cdot\mathbb{E} \bigl[ \exp( p z \alpha_{N-1}^N ) | \mathcal{F}_{ {
(N-1) T }/{ N } } \bigr] \Biggr]
\nonumber\hspace*{-25pt}\\
&\leq&
\mathbb{E}\Biggl[ \exp\Biggl( p z \sum_{ k = 0 }^{ N - 2 } \alpha_k^N \Biggr)
\Biggr] \cdot\exp\biggl( \frac{ p^2 T ( c + \| \sigma(0) \| )^2 }{ N } \biggr)\hspace*{-25pt}
\\
&\leq&
\cdots\leq \exp\bigl( p^2 T \bigl( c + \| \sigma(0) \| \bigr)^2
\bigr)\nonumber\hspace*{-25pt}
\end{eqnarray}
for all $ N \in\mathbb{N} $, $ p \in[1, \infty) $ and all $ z \in\{ -1,
1 \} $. Combining~(\ref{eqsupproofA}) and~(\ref{eqsupproofC}) then
gives that
\begin{eqnarray*}
&&
\sup_{ z \in\{ -1, 1 \} } \sup_{ N \in\mathbb{N} } \Biggl\| \sup_{ n \in\{ 0,
1, \ldots, N \} } \exp\Biggl( z \sum_{ k = 0 }^{ n - 1 } \alpha_k^N \Biggr)
\Biggr\|_{
L^p( \Omega; \mathbb{R} ) } \\
&&\qquad\leq2 \exp\bigl( p^2 T \bigl( c + \| \sigma(0) \|
\bigr)^2 \bigr) < \infty
\end{eqnarray*}
for all $ p \in[2, \infty) $ and this
completes the proof of Lemma~\ref{lemsupsupsup}.

%
\subsection{\texorpdfstring{Proof of Lemma \protect\ref{lembdmomodom}}{Proof of Lemma 3.5}}

First of all, H{\"o}lder's inequality shows that
%
%
\begin{eqnarray}
\label{eqbdmomproofA}
&&
\mathop{\sup_{N \in\mathbb{N},}}_{N \geq8 \lambda
p T} \Bigl\| \sup_{ n \in\{ 0, 1, \ldots, N \} } D_n^N \Bigr\|_{ L^p( \Omega;
\mathbb{R} ) }
\nonumber\hspace*{-35pt}\\
&&\qquad\leq
e^{ \lambda} \bigl( \lambda+ \| \xi\|_{ L^{ 4 p }( \Omega; \mathbb{R}^d )
} \bigr) \Biggl( \mathop{\sup_{N \in\mathbb{N},}}_{N \geq8 \lambda p T}
\Biggl\| \exp\Biggl(
\lambda\sum_{ k = 0 }^{ N - 1 } \| \Delta W_k^N \|^2 \Biggr) \Biggr\|_{ L^{ 2 p }(
\Omega; \mathbb{R} ) } \Biggr)\hspace*{-35pt}
\\
&&\qquad\quad{}\times
\Biggl( \sup_{ N \in\mathbb{N} } \Biggl\| \sup_{ n \in\{ 0, 1, \ldots, N \} } \exp
\Biggl( \sup_{ u \in\{ 0, 1, \ldots, n \} } \sum_{ k = u }^{ n - 1 }
\alpha_k^N \Biggr) \Biggr\|_{ L^{ 4 p }( \Omega; \mathbb{R} ) } \Biggr) \nonumber\hspace*{-35pt}
\end{eqnarray}
for all $ p \in[ 1, \infty) $. Moreover, again H{\"o}lder's inequality
implies that
%
%
\begin{eqnarray}
\label{eqbdmomproofB}
&&
\Biggl\| \sup_{ n \in\{ 0, 1, \ldots, N \} } \exp\Biggl(
\sup_{ u \in\{ 0, 1, \ldots, n \} } \sum_{ k = u }^{ n - 1 } \alpha_k^N
\Biggr) \Biggr\|_{ L^{ 4 p }( \Omega; \mathbb{R} ) }
\nonumber\\
&&\qquad\leq
\Biggl\| \sup_{ n \in\{ 0, 1, \ldots, N \} } \exp\Biggl( \sum_{ k = 0
}^{ n - 1 } \alpha_k^N \Biggr) \Biggr\|_{ L^{ 8p }( \Omega; \mathbb{R} ) }
\\
&&\qquad\quad{}\times\Biggl\|
\sup_{ u \in\{ 0, 1, \ldots, N \} } \exp\Biggl( -\sum_{ k = 0 }^{ u - 1 }
\alpha_k^N \Biggr) \Biggr\|_{ L^{ 8p }( \Omega; \mathbb{R} ) }\nonumber
\end{eqnarray}
for all $ N \in\mathbb{N} $ and all $ p \in[1, \infty) $. Putting
(\ref{eqbdmomproofB}) into~(\ref{eqbdmomproofA}) then gives that
%
%
\begin{eqnarray}
\label{eqbdmomproof}
&&\mathop{\sup_{N \in\mathbb{N},}}_{N \geq8 \lambda
p T} \Bigl\| \sup_{ n \in\{ 0, 1, \ldots, N \} } D_n^N \Bigr\|_{ L^p( \Omega;
\mathbb{R}^d ) } \nonumber\hspace*{-30pt}\\[-3pt]
&&\qquad\leq e^{ \lambda} \bigl( \lambda+ \| \xi\|_{
L^{ 4 p }( \Omega; \mathbb{R} ) } \bigr) \Biggl( \mathop{\sup_{N
\in\mathbb{N},}}_{N \geq8 \lambda p T} \Biggl\| \exp\Biggl( \lambda\sum_{ k = 0
}^{ N - 1 } \| \Delta W_k^N \|^2 \Biggr) \Biggr\|_{ L^{ 2 p }( \Omega; \mathbb{R} )
} \Biggr)\hspace*{-30pt}
\\[-3pt]
&&\qquad\quad{}\times\Biggl( \sup_{ z \in\{-1,1\} }
\sup_{ N \in\mathbb{N} } \Biggl\| \sup_{ n
\in\{ 0, 1, \ldots, N \} } \exp\Biggl( z \sum_{ k = 0 }^{ n - 1 } \alpha_k^N
\Biggr) \Biggr\|_{ L^{ 8p }( \Omega; \mathbb{R} ) } \Biggr)^{ 2 } \nonumber\hspace*{-30pt}
\end{eqnarray}
for all $ p \in[1, \infty) $. Combining Lemmas~\ref{lemexplemma},
\ref{lemsupsupsup} and~(\ref{eqbdmomproof}) finally completes the proof
of Lemma~\ref{lembdmomodom}.

%
\subsection{\texorpdfstring{Proof of Lemma \protect\ref{lOmegaN}}{Proof of Lemma 3.6}}
The subadditivity of the probability measure~$ \mathbb{P}[ \cdot] $ and
Markov's inequality yield that
%
%
\begin{eqnarray}\label{eqOmegaN}
\mathbb{P}[ ( \Omega_N^N )^c ]
&\leq& \mathbb{P}\Bigl[ \sup_{ n \in\{ 0, 1,
\ldots, N-1 \} } D_n^N
>
N^{ { 1 }/({ 2c }) } \Bigr] + N \cdot\mathbb{P}[ \| W_{ { T }/{ N } }
\|
> 1
] \nonumber\hspace*{-35pt}\\[-3pt]
&\leq& \mathbb{E}\Bigl[ \sup_{ n \in\{ 0, 1, \ldots, N-1 \} } |
D_n^N |^q \Bigr] \cdot N^{ { -q }/({ 2c }) } + N \cdot\mathbb{P}\bigl[ \| W_{
T } \|
> \sqrt{N}
\bigr]\hspace*{-35pt} \\[-3pt]
&\leq& \mathbb{E}\Bigl[ \sup_{ n \in\{ 0, 1, \ldots, N \} } |
D_n^N |^q \Bigr] \cdot N^{ { -q }/({ 2c }) } + \mathbb{E}[ \| W_{ T } \|^q
] \cdot N^{ ( 1 - { q }/{ 2 } ) }\nonumber\hspace*{-35pt}
\end{eqnarray}
for all $ N \in\mathbb{N} $ and all $ q \in[1, \infty) $. Combining
Lemma~\ref{lembdmomodom} and~(\ref{eqOmegaN}) then shows~(\ref{eqOmegaNgoal}). This completes the proof of Lemma~\ref{lOmegaN}.

\subsection{\texorpdfstring{Proof of Lemma \protect\ref{lBDG}}{Proof of Lemma 3.7}}
\label{seclBDG}
Lemma~\ref{lBDG} is an immediate consequence of Doob's maximal
inequality, of the Burkholder--Davis--Gundy type inequality in
Lem\-ma~7.7 of Da Prato and Zabczyk~\cite{dz92} and of the triangle inequality.
For completeness we now present the proof of Lemma~\ref{lBDG}.
\begin{pf*}{Proof of Lemma~\ref{lBDG}}
Doob's maximal inequality (see, e.g., Da Prato and Zabczyk~\cite{dz92},
Theorem 3.8), Da Prato and Zabczyk~\cite{dz92}, Lemma 7.7, and the
triangle inequality give
%
%
\begin{eqnarray}
\label{eqBDG2}
&&\E\biggl[\sup_{s\in[0,t]} \biggl\| \int_0^s Z_u \,dW_u \biggr\|^p \biggr] \nonumber\\[-3pt]
&&\qquad \leq
\biggl(\frac{ p }{ p - 1 } \biggr)^{ p } \biggl( \frac{ p ( p - 1 ) }{2} \biggr)^{ {p}/{2}}
\Biggl( \int_0^t \Biggl\| \sum_{i=1}^m\| Z_s \vec{ e }_i \|^2 \Biggr\|_{ L^{ { p }/{
2 } }( \Omega; \mathbb{R} ) } \,ds \Biggr)^{ {p}/{2} } \nonumber\\[-3pt]
&&\qquad \leq
\biggl( \frac{ p }{ p - 1 } \biggr)^{ p } \biggl( \frac{ p ( p - 1 ) }{2} \biggr)^{
{p}/{2}} \Biggl( \int_0^t \sum_{i=1}^m \| \| Z_s \vec{ e }_i \|^2 \|_{
L^{ { p }/{ 2 } }( \Omega; \mathbb{R} ) } \,ds \Biggr)^{ {p}/{2} }
\\[-3pt]
&&\qquad = | p |^p \biggl( \frac{ p }{ 2 ( p - 1 ) } \biggr)^{ {p}/{2}} \Biggl(
\int_0^t \sum_{i=1}^m \| Z_s \vec{ e }_i \|^2_{ L^{ p }( \Omega;
\mathbb{R}^k ) } \,ds \Biggr)^{ {p}/{2} } \nonumber\\[-3pt]
&&\qquad \leq | p |^p \Biggl(
\int_0^t \sum_{i=1}^m \| Z_s \vec{ e }_i \|^2_{ L^{ p }( \Omega;
\mathbb{R}^k ) } \,ds \Biggr)^{ {p}/{2} }\nonumber
\end{eqnarray}
for all $ t \in[0,T] $ and all $ p \in[2,\infty) $. This completes the
proof of Lemma~\ref{lBDG}.\vspace*{-3pt}
\end{pf*}

\subsection{\texorpdfstring{Proof of Lemma \protect\ref{cBDG}}{Proof of Lemma 3.8}}
\label{seccBDG}

Let $ \bar{Z}^N \dvtx[0,T] \times\Omega\rightarrow\mathbb{R}^{ k
\times m } $, $ N \in\mathbb{N} $, be a sequence of stochastic
processes defined by $ \bar{Z}^N_s := Z_l^N $ for all
$s\in[\frac{lT}{N},\frac{(l+1)T}{N})$, $l\in\{0,1,\ldots,N-1\}$ and
all $N\in\N$. The Burkholder--Davis--Gundy type inequality in Lemma
\ref{lBDG} then shows
\begin{eqnarray*}
\hspace*{-4pt}&& \Biggl\| \sup_{j\in\{0,1,\ldots,n\}} \Biggl\| \sum_{l=0}^{j-1}Z_l^N\Delta W_l^N
\Biggr\| \Biggr\|_{L^p(\Omega;\R)}
\\[-3pt]
\hspace*{-4pt}&&\quad=
\biggl\| \sup_{j\in\{0,1,\ldots,n\}} \biggl\| \int_0^{ { j T }/{ N } }
\bar{Z}_u
\,dW_u \biggr\| \biggr\|_{L^p(\Omega;\R)}
\leq\biggl\| \sup_{ s \in[0,{ n T }/{ N }
] } \biggl\| \int_0^{ s } \bar{Z}_u
\,dW_u \biggr\| \biggr\|_{L^p(\Omega;\R)}
\\[-3pt]
\hspace*{-4pt}&&\quad\leq
p \Biggl( \int_0^{ { n T }/{ N } } \sum_{ i = 1 }^m \| \bar{Z}_s \vec{ e
}_i \|^2_{ L^p( \Omega; \mathbb{R}^k ) } \,ds \Biggr)^{ { 1 }/{ 2 } } =
p\Biggl(\sum_{l=0}^{n-1}\sum_{i=1}^m \|Z_l^N \vec{e}_i
\|_{L^p(\Omega;\R^{k})}^2 \frac{T}{N} \Biggr)^{{1}/{2}}
\end{eqnarray*}
for all $ n \in\{ 0, 1, \ldots, N \} $,
$ N \in\mathbb{N} $
and all $ p \in[2,\infty) $.
This completes the proof
of Lemma~\ref{cBDG}.\vspace*{-3pt}

\subsection{\texorpdfstring{Proof of Lemma \protect\ref{cmomEuler}}{Proof of Lemma 3.9}}

In order to show Lemma~\ref{cmomEuler} we first represent the numerical
approximations~(\ref{eqscheme}) in an appropriate way. More formally,
we have that
%
%
\begin{eqnarray}
\label{eqrep} Y_n^N &=& \xi+ \sum_{ k = 0 }^{ n - 1 } \frac{
({ T}/{ N }) \mu( Y_k^N ) }{ 1 + ({ T}/{ N })\| \mu( Y_k^N ) \| } + \sum_{
k = 0 }^{ n - 1 } \sigma( Y_k^N ) \Delta W_k^N \nonumber\\[-3pt]
&=& \xi+
\sigma( 0 ) W_{ {n T}/{N} } + \sum_{ k = 0 }^{ n - 1 } \frac{
({ T}/{ N }) \mu( Y_k^N ) }{ 1 + ({ T}/{ N })\| \mu( Y_k^N ) \| }\\[-3pt]
&&{}+ \sum_{ k = 0 }^{ n - 1 } \bigl( \sigma( Y_k^N ) - \sigma( 0 ) \bigr) \Delta
W_k^N\nonumber
\end{eqnarray}
$ \mathbb{P} $-a.s. for all $ n \in\{ 0,1, \ldots,N \} $ and all $ N
\in\mathbb{N} $. The Burkholder--Davis--Gundy type
inequality in Lemma~\ref{cBDG}
then gives that
%
%
\begin{eqnarray}
&& \| Y_n^N \|_{ L^p( \Omega; \mathbb{R}^d ) } \nonumber\\
&&\qquad\leq\|
\xi\|_{ L^p( \Omega; \mathbb{R}^d ) } + \|
\sigma( 0 )
W_{ { n T }/{ N } } \|_{ L^p( \Omega; \mathbb{R}^d ) }\nonumber\\
&&\qquad\quad{} + \sum_{ k =
0 }^{ n - 1 } \biggl\| \frac{ ({ T}/{ N }) \mu( Y_k^N ) }{ 1 + ({ T}/{ N }) \|
\mu( Y_k^N ) \| } \biggr\|_{ L^p( \Omega; \mathbb{R}^d ) }
\nonumber\\[-8pt]\\[-8pt]
&&\qquad\quad{}+ \Biggl\|
\sum_{ k = 0 }^{ n - 1 }
\bigl( \sigma( Y_k^N ) - \sigma( 0 ) \bigr)
\Delta W_k^N \Biggr\|_{ L^p( \Omega; \mathbb{R}^d ) } \nonumber\\
&&\qquad\leq\|
\xi\|_{ L^p( \Omega; \mathbb{R}^d ) } + p\Biggl( \frac{ nT }{ N } \sum_{ i =
1 }^{ m } \| \sigma_i( 0 ) \|^2 \Biggr)^{ { 1 }/{ 2 } } + N \nonumber\\
&&\qquad\quad{}
+ p\Biggl( \sum_{ k = 0 }^{ n - 1 } \sum_{ i = 1 }^{ m } \| \sigma_i( Y_k^N )
- \sigma_i( 0 ) \|^2_{ L^p( \Omega; \mathbb{R}^d ) } \frac{ T }{ N }
\Biggr)^{ { 1 }/{ 2 } }\nonumber
\end{eqnarray}
and the global Lipschitz continuity of $ \sigma$ therefore shows that
\begin{eqnarray*}
\| Y_n^N \|_{ L^p( \Omega; \mathbb{R}^d ) }^2
&\leq&
2\bigl( \| \xi\|_{ L^p( \Omega; \mathbb{R}^d ) } + p \sqrt{ T m } \|
\sigma( 0 ) \| + N \bigr)^2\\
&&{} + \frac{ 2 p^2 m T c^2 }{ N } \Biggl( \sum_{ k = 0 }^{
n - 1 } \| Y_k^N \|^2_{ L^p( \Omega; \mathbb{R}^d ) } \Biggr)
\end{eqnarray*}
for all $ n \in\{ 0, 1, \ldots, N \} $, $ N \in\mathbb{N} $ and all $ p
\in[2, \infty) $. In the next step Gronwall's lemma gives that
%
%
\begin{eqnarray}
\label{eqboundY1}
&&
\sup_{ n \in\{ 0,1, \ldots, N \} } \| Y_n^N \|_{ L^p(
\Omega; \mathbb{R}^d ) } \nonumber\\[-8pt]\\[-8pt]
&&\qquad\leq\sqrt{ 2 } e^{ p^2 m T c^2 } \bigl( \| \xi\|_{
L^p( \Omega; \mathbb{R}^d ) } + p \sqrt{ T m } \| \sigma( 0 ) \| + N
\bigr)\nonumber
\end{eqnarray}
for all $ N \in\mathbb{N} $ and all $ p \in[2,\infty) $. Of course,
(\ref{eqboundY1}) does not prove Lemma~\ref{cmomEuler} due to the $N
\in\mathbb{N}$ on the right-hand side of~(\ref{eqboundY1}). However,
exploiting~(\ref{eqboundY1}) in an appropriate bootstrap argument will
enable us to establish Lemma~\ref{cmomEuler}. More formally,
H{\"o}lder's inequality, estimate~(\ref{eqboundY1}) and Lemma
\ref{lOmegaN} show that
%
%
\begin{eqnarray}
\label{eqboundY2}
&&
\sup_{ N \in\mathbb{N} } \sup_{ n \in\{ 0,
1, \ldots, N \} } \bigl\| \mathbh{1}_{ ( \Omega_n^N )^c } Y_n^N \bigr\|_{ L^p(
\Omega; \mathbb{R}^d ) }
\nonumber\\
&&\qquad\leq
\sup_{ N \in\mathbb{N} } \sup_{ n \in\{ 0, 1, \ldots, N \} }
\bigl( \bigl\| \mathbh{1}_{ ( \Omega_N^N )^c } \bigr\|_{ L^{2p}
( \Omega; \mathbb{R}
) } \| Y_n^N \|_{ L^{2p}( \Omega; \mathbb{R}^d ) } \bigr)
\nonumber\\
&&\qquad\leq
\Bigl( \sup_{ N \in\mathbb{N} } \bigl( N \cdot\bigl\| \mathbh{1}_{ (
\Omega_N^N )^c } \bigr\|_{ L^{2p}( \Omega; \mathbb{R} ) } \bigr) \Bigr)
\nonumber\\
&&\qquad\quad{}\times\Bigl( \sup_{ N
\in\mathbb{N} } \sup_{ n \in\{ 0, 1, \ldots, N \} } \bigl( N^{ -1 } \cdot\|
Y_n^N \|_{ L^{2p}( \Omega; \mathbb{R}^d ) } \bigr) \Bigr)
\\
&&\qquad\leq
\sqrt{2} e^{ 4 p^2 m T c^2 } \Bigl( \sup_{ N \in\mathbb{N} } N^{
2p } \cdot\mathbb{P}[ ( \Omega_N^N )^c ] \Bigr)^{ { 1 }/({ 2 p }) } \nonumber\\
&&\qquad\quad{}\times\Bigl( \|
\xi\|_{ L^{ 2p }( \Omega; \mathbb{R}^d ) } + 2 p \sqrt{ T m } \|
\sigma( 0 ) \| + 1 \Bigr)
\nonumber\\
&&\qquad
< \infty\nonumber
\end{eqnarray}
for all $ p \in[2, \infty) $. Additionally,
Lemmas~\ref{ldominator} and~\ref{lembdmomodom} give that
%
%
\begin{eqnarray}
\label{eqboundY3}
&&
\mathop{\sup_{N \in\mathbb{N},}}_{N \geq8 \lambda p T
} \sup_{ n \in\{ 0, 1, \ldots, N \} } \| \mathbh{1}_{ \Omega_n^N }
Y_n^N \|_{ L^p( \Omega; \mathbb{R}^d ) } \nonumber\\[-8pt]\\[-8pt]
&&\qquad\leq\mathop{\sup_{N
\in\mathbb{N},}}_{N \geq8 \lambda p T} \sup_{ n \in\{ 0, 1, \ldots, N
\} } \| D_n^N \|_{ L^p( \Omega; \mathbb{R} ) } < \infty\nonumber
\end{eqnarray}
for all $ p \in[1, \infty) $. Combining~(\ref{eqboundY2}) and
(\ref{eqboundY3}) finally completes the proof of Lem\-ma~\ref{cmomEuler}.

\subsection{\texorpdfstring{Proof of Lemma \protect\ref{lemmubound}}{Proof of Lemma 3.10}}
\label{secproofmubound}

Lemma~\ref{lemmubound}
follows from Lemma~\ref{cmomEuler}. More precisely, the estimate $ \|
\mu(x) \| \leq( 2 c + \| \mu(0) \| ) ( 1 + \| x \|^{ (c+1) } ) $ for
all $ x \in\mathbb{R}^d $ and Lemma~\ref{cmomEuler} give that
%
%
\begin{eqnarray}
&&
\sup_{ N \in\mathbb{N} } \sup_{ n \in\{ 0, 1, \ldots, N \} } \| \mu(
Y^N_n ) \|_{ L^p( \Omega; \mathbb{R}^d )}
\nonumber\\
&&\qquad\leq\bigl( 2 c + \| \mu(0) \| \bigr) \Bigl( 1 + \sup_{ N \in\mathbb{N} } \sup
_{ n
\in\{ 0, 1, \ldots, N \} } \| Y^N_n \|^{ (c + 1) }_{ L^{ p (c+1) }(
\Omega; \mathbb{R}^d ) } \Bigr) \\
&&\qquad< \infty\nonumber
\end{eqnarray}
for all $ p \in[1,\infty) $. Additionally, the inequality $ \|
\sigma(x) \| \leq c \| x \| + \| \sigma( 0 ) \| $ for all $ x
\in\mathbb{R}^d $ and again Lemma~\ref{cmomEuler} show that
%
%
\begin{eqnarray}
&&\sup_{ N \in\mathbb{N} } \sup_{ n \in\{ 0, 1, \ldots, N \} } \| \sigma(
Y^N_n ) \|_{ L^p( \Omega; \mathbb{R}^d ) }
\nonumber\\
&&\qquad\leq c \Bigl( \sup_{ N \in\mathbb{N} } \sup_{ n \in\{ 0, 1, \ldots,
N \} }
\| Y^N_n \|_{ L^{ p }( \Omega; \mathbb{R}^d ) } \Bigr) + \| \sigma(0) \|\\
&&\qquad <
\infty\nonumber
\end{eqnarray}
for all $ p \in[1,\infty) $. This completes the proof of Lemma
\ref{lemmubound}.

\subsection{\texorpdfstring{Proof of Theorem \protect\ref{thmmain}}{Proof of Theorem 3.11}}
\label{secproofthm}

Using Lemmas~\ref{cmomEuler} and~\ref{lemmubound} we now establish
inequality~(\ref{eqmainresult}). To this end we use the notation
\[
\lfloor t \rfloor_N := \max\biggl\{ s \in\biggl\{ 0, \frac{ T }{ N } , \frac{ 2 T
}{ N } , \ldots, \frac{ (N-1) T }{ N }, T \biggr\} \dvtx s \leq t \biggr\}
\]
for all $ t \in[0,T] $ and all $ N \in\mathbb{N} $. In this notation,
equation~(\ref{eqtimecont}) reads as
%
%
\begin{eqnarray}
\label{eqtildeapprox} \bar{Y}_s^N &=& \xi+ \int_0^s \frac{ \mu(
\bar{Y}_{ \lfloor u \rfloor_N }^N ) }{ 1 + ({ T }/{ N }) \| \mu(
\bar{Y}_{ \lfloor u \rfloor_N }^N ) \| } \,du \nonumber\\[-8pt]\\[-8pt]
&&{}+ \int_0^s \sigma\bigl(
\bar{Y}_{ \lfloor u \rfloor_N }^N \bigr) \,dW_u\nonumber
\end{eqnarray}
for all $ s \in[0,T] $ $ \mathbb{P} $-a.s. and all $ N \in\mathbb{N} $.
Our goal is then to estimate the quantity $ \mathbb{E}[ \sup_{ s
\in[0,t] } \| X_s - \bar{Y}_s^N \|^p ] $ for $ t \in[0,T] $ and $ p
\in[1,\infty) $. To this end note that~(\ref{eqSDEsol}) and
(\ref{eqtildeapprox}) imply that
\begin{eqnarray*}
X_s - \bar{Y}_s^N
&=&
\int_0^s \biggl( \mu( X_u ) - \frac{ \mu( \bar{Y}_{ \lfloor u \rfloor_N }^N )
}{ 1 + ({ T }/{ N }) \| \mu( \bar{Y}_{ \lfloor u \rfloor_N }^N ) \| }
\biggr) \,du\\
&&{} + \sum_{ i = 1 }^m \int_0^s \bigl( \sigma_i( X_u ) - \sigma_i\bigl(
\bar{Y}_{ \lfloor u \rfloor_N }^N \bigr) \bigr) \,dW_u^{(i)}
\end{eqnarray*}
for all $ s \in[0,T] $ $ \mathbb{P} $-a.s. and all $ N \in\mathbb{N} $.
It{\^o}'s formula hence gives that
\begin{eqnarray*}
\| X_s - \bar{Y}_s^N \|^2 &=& 2 \int_0^s \langle X_u - \bar{Y}_u^N, \mu( X_u )
- \mu( \bar{Y}_u^N )
\rangle \,du
\\
&&{} +
2 \int_0^s \bigl\langle X_u - \bar{Y}_u^N, \mu( \bar{Y}_u^N ) - \mu\bigl( \bar{Y}_{
\lfloor u \rfloor_N }^N \bigr)
\bigr\rangle \,du
\\
&&{} +
\frac{ 2 T }{ N } \int_0^s \biggl\langle X_u - \bar{Y}_u^N, \frac{ \mu( \bar{Y}_{
\lfloor u \rfloor_N }^N ) \| \mu( \bar{Y}_{ \lfloor u \rfloor_N }^N )
\| }{ 1 + ({T}/{N}) \| \mu( \bar{Y}_{ \lfloor u \rfloor_N }^N ) \| }
\biggr\rangle \,du
\\
&&{} +
2 \sum_{ i = 1 }^m \int_0^s \bigl\langle X_u - \bar{Y}_u^N, \sigma_i( X_u ) -
\sigma_i\bigl( \bar{Y}_{ \lfloor u \rfloor_N }^N \bigr)
\bigr\rangle \,dW^{(i)}_u
\\
&&{} +
\sum_{ i = 1 }^m \int_0^s \bigl\| \sigma_i( X_u ) - \sigma_i\bigl( \bar{Y}_{
\lfloor u \rfloor_N }^N \bigr) \bigr\|^2 \,du
\end{eqnarray*}
and the inequality $ ( a + b )^2 \leq2a^2 + 2b^2 $ for all $ a, b
\in\mathbb{R} $, the estimate $ a \leq| a | $ for all $ a \in\mathbb{R}
$ and the Cauchy--Schwarz inequality therefore yield that
\begin{eqnarray*}
\| X_s - \bar{Y}_s^N \|^2 &\leq& ( 2 c + 2c^2 m ) \int_0^s \| X_u -
\bar{Y}_u^N \|^2 \,du
\\
&&{} +
2 \int_0^s \| X_u - \bar{Y}_u^N \| \bigl\| \mu( \bar{Y}_u^N ) - \mu\bigl(
\bar{Y}_{ \lfloor u \rfloor_N }^N \bigr) \bigr\| \,du
\\
&&{} +
\frac{ 2 T }{ N } \int_0^s \| X_u - \bar{Y}_u^N \| \bigl\| \mu\bigl( \bar{Y}_{
\lfloor u \rfloor_N }^N \bigr) \bigr\|^2 \,du
\\
&&{} +
2 \Biggl| \sum_{ i = 1 }^m \int_0^s \bigl\langle X_u - \bar{Y}_u^N, \sigma_i( X_u ) -
\sigma_i\bigl( \bar{Y}_{ \lfloor u \rfloor_N }^N \bigr)
\bigr\rangle \,dW^{(i)}_u
\Biggr|
\\
&&{} +
2 c^2 m \int_0^s \bigl\| \bar{Y}_u^N - \bar{Y}_{ \lfloor u \rfloor_N }^N
\bigr\|^2 \,du
\end{eqnarray*}
for all $ s \in[0,T] $ $ \mathbb{P} $-a.s. and all $ N \in\mathbb{N} $.
The inequality $ a \cdot b \leq\frac{a^2}{2} + \frac{b^2}{2} $ for all
$ a, b \in\mathbb{R} $ then shows that
\begin{eqnarray*}
&&{\sup_{ s \in[0,t] }} \| X_s - \bar{Y}_s^N \|^2
\\
&&\qquad\leq
2 ( c + c^2 m + 1 ) \int_0^t \| X_s - \bar{Y}_s^N \|^2 \,ds + \int_0^T \bigl\|
\mu( \bar{Y}_s^N ) - \mu\bigl( \bar{Y}_{ \lfloor s \rfloor_N }^N \bigr) \bigr\|^2 \,ds
\\
&&\qquad\quad{} +
\frac{ T^2 }{ N^2 } \int_0^T \bigl\| \mu\bigl( \bar{Y}_{ \lfloor s \rfloor_N }^N
\bigr) \bigr\|^4 \,ds + 2 c^2 m \int_0^T \bigl\| \bar{Y}_s^N - \bar{Y}_{ \lfloor s
\rfloor_N }^N \bigr\|^2 \,ds
\\
&&\qquad\quad{} +
2 \sup_{ s \in[0,t] } \Biggl| \sum_{ i = 1 }^m \int_0^s \bigl\langle X_u - \bar{Y}_u^N,
\sigma_i( X_u ) - \sigma_i\bigl( \bar{Y}_{ \lfloor u \rfloor_N }^N \bigr)
\bigr\rangle \,dW^{(i)}_u
\Biggr|
\end{eqnarray*}
$ \mathbb{P} $-a.s. for
all $ t \in[0,T] $
and all
$ N \in\mathbb{N} $. The
Burkholder--Davis--Gundy type inequality in Lemma~\ref{lBDG} hence
yields that
%
%
\begin{eqnarray}
\label{eqafterBDG}
&&\Bigl\| \sup_{ s \in[0,t] } \| X_s - \bar{Y}_s^N \|^2
\Bigr\|_{ L^{ {p}/{2} }( \Omega; \mathbb{R} ) } \nonumber\\
&&\qquad\leq
2 ( c +
c^2 m + 1 ) \int_0^t \| X_s - \bar{Y}_s^N \|^2_{ L^p( \Omega;
\mathbb{R}^d ) } \,ds \nonumber\\
&&\qquad\quad{} + \int_0^T \bigl\| \mu( \bar{Y}_s^N ) -
\mu\bigl( \bar{Y}_{ \lfloor s \rfloor_N }^N \bigr) \bigr\|^2_{ L^p( \Omega;
\mathbb{R}^d ) } \,ds \\
&&\qquad\quad{} + \frac{ T^2 }{ N^2 } \int_0^T \bigl\| \mu\bigl(
\bar{Y}_{ \lfloor s \rfloor_N }^N \bigr) \bigr\|^4_{ L^{2p}( \Omega; \mathbb{R}^d
) } \,ds + 2 c^2 m \int_0^T \bigl\| \bar{Y}_s^N - \bar{Y}_{ \lfloor s
\rfloor_N }^N \bigr\|^2_{ L^p( \Omega; \mathbb{R}^d ) } \,ds \nonumber\\
&&\qquad\quad{} + p
\Biggl( \sum_{ i = 1 }^m \int_0^t \bigl\| \bigl\langle X_s - \bar{Y}_s^N,
\sigma_i( X_s ) -
\sigma_i\bigl( \bar{Y}_{ \lfloor s \rfloor_N }^N \bigr)
\bigr\rangle
\bigr\|^2_{ L^{{p}/{2}}( \Omega; \mathbb{R} ) } \,ds \Biggr)^{ {1}/{2}
}\nonumber
\end{eqnarray}
%
for all $ t \in[0,T]$, $ N \in\mathbb{N} $ and all $ p \in[4,
\infty
) $.
Next the Cauchy--Schwarz inequality, the H{\"o}lder inequality
and again the inequality $ a \cdot b \leq\frac{a^2}{2} + \frac{b^2}{2}
$ for all $ a, b \in\mathbb{R} $
imply that
%
%
\begin{eqnarray}
\label{eq2afterBDG}\quad
&&p \Biggl( \sum_{ i = 1 }^m \int_0^t \bigl\| \bigl\langle X_s -
\bar{Y}_s^N, \sigma_i( X_s ) - \sigma_i\bigl( \bar{Y}_{ \lfloor s \rfloor_N
}^N \bigr)
\bigr\rangle
\bigr\|^2_{ L^{{p}/{2}}( \Omega; \mathbb{R} ) } \,ds \Biggr)^{ {1}/{2} }
\nonumber\\
&&\qquad\leq p \Biggl( \sum_{ i = 1 }^m \int_0^t \| X_s - \bar{Y}_s^N
\|_{ L^p( \Omega; \mathbb{R}^d ) }^2 \bigl\| \sigma_i( X_s ) - \sigma_i\bigl(
\bar{Y}_{ \lfloor s \rfloor_N }^N \bigr) \bigr\|^2_{ L^p( \Omega; \mathbb{R}^d )
} \,ds \Biggr)^{ {1}/{2} } \nonumber\\
&&\qquad\leq p \Bigl( {\sup_{ s \in[0,t] }} \|
X_s - \bar{Y}_s^N \|_{ L^p( \Omega; \mathbb{R}^d ) } \Bigr) \biggl( c^2 m \int_0^t
\bigl\| X_s - \bar{Y}_{ \lfloor s \rfloor_N }^N \bigr\|^2_{ L^p( \Omega;
\mathbb{R}^d ) } \,ds \biggr)^{ {1}/{2} } \\
&&\qquad\leq{\frac{1}{2}
\sup_{ s \in[0,t] }} \| X_s - \bar{Y}_s^N \|^2_{ L^p( \Omega;
\mathbb{R}^d ) } + \frac{ p^2 c^2 m }{ 2 } \int_0^t \bigl\| X_s - \bar{Y}_{
\lfloor s \rfloor_N }^N \bigr\|^2_{ L^p( \Omega; \mathbb{R}^d ) } \,ds
\nonumber\\[-2pt]
&&\qquad\leq\frac{1}{2} \Bigl\| \sup_{ s \in[0,t] } \| X_s -
\bar{Y}_s^N \|\Bigr\|^2_{ L^p( \Omega; \mathbb{R} ) } + \frac{ p^2 c^2 m }{
2 } \int_0^t \bigl\| X_s - \bar{Y}_{ \lfloor s \rfloor_N }^N \bigr\|^2_{ L^p(
\Omega; \mathbb{R}^d ) } \,ds\nonumber
\end{eqnarray}
for all $ t\!\in\![0,T]$, $ N\!\in\!\mathbb{N} $ and all $ p\!\in\![4, \infty)
$. Inserting inequality~(\ref{eq2afterBDG}) into~(\ref{eqafterBDG})~and
applying the estimate $ ( a + b )^2 \leq2a^2 + 2b^2 $ for all $ a, b
\in\mathbb{R} $ then yields that
\begin{eqnarray*}
&&\Bigl\| \sup_{ s \in[0,t] } \| X_s - \bar{Y}_s^N \| \Bigr\|_{ L^{p }( \Omega;
\mathbb{R} ) }^2
\\[-2pt]
&&\qquad\leq
2 \biggl( c + c^2 m + 1 + \frac{ p^2 c^2 m }{ 2 } \biggr) \int_0^t \| X_s -
\bar{Y}_s^N \|^2_{ L^p( \Omega; \mathbb{R}^d ) } \,ds
\\[-2pt]
&&\qquad\quad{} +
\int_0^T \bigl\| \mu( \bar{Y}_s^N ) - \mu\bigl( \bar{Y}_{ \lfloor s \rfloor_N }^N
\bigr) \bigr\|^2_{ L^p( \Omega; \mathbb{R}^d ) } \,ds + \frac{ T^2 }{ N^2 }
\int_0^T \bigl\| \mu\bigl( \bar{Y}_{ \lfloor s \rfloor_N }^N \bigr) \bigr\|^4_{ L^{2p}(
\Omega; \mathbb{R}^d ) } \,ds
\\[-2pt]
&&\qquad\quad{} +
( 2 c^2 m + p^2 c^2 m ) \int_0^T \bigl\| \bar{Y}_s^N - \bar{Y}_{ \lfloor s
\rfloor_N }^N \bigr\|^2_{ L^p( \Omega; \mathbb{R}^d ) } \,ds
\\[-2pt]
&&\qquad\quad{} +
\frac{1}{2} \Bigl\| \sup_{ s \in[0,t] } \| X_s - \bar{Y}_s^N \| \Bigr\|^2_{ L^p(
\Omega; \mathbb{R} ) }
\end{eqnarray*}
%
and therefore, we obtain
that
\begin{eqnarray*}
\hspace*{-3pt}&&\frac{1}{2} \Bigl\| \sup_{ s \in[0,t] } \| X_s - \bar{Y}_s^N \| \Bigr\|^2_{ L^p(
\Omega; \mathbb{R} ) }
\\[-2pt]
\hspace*{-3pt}&&\qquad\leq
2 ( c + 1 + p^2 c^2 m ) \int_0^t \| X_s - \bar{Y}_s^N \|^2_{ L^p(
\Omega; \mathbb{R}^d ) } \,ds
\\[-2pt]
\hspace*{-3pt}&&\qquad\quad{} +
\int_0^T \bigl\| \mu( \bar{Y}_s^N ) - \mu\bigl( \bar{Y}_{ \lfloor s \rfloor_N }^N
\bigr) \bigr\|^2_{ L^p( \Omega; \mathbb{R}^d ) } \,ds + \frac{ T^2 }{ N^2 }
\int_0^T \bigl\| \mu\bigl( \bar{Y}_{ \lfloor s \rfloor_N }^N \bigr) \bigr\|^4_{ L^{2p}(
\Omega; \mathbb{R}^d ) } \,ds
\\[-2pt]
\hspace*{-3pt}&&\qquad\quad{} +
( 2 c^2 m + p^2 c^2 m ) \int_0^T \bigl\| \bar{Y}_s^N - \bar{Y}_{ \lfloor s
\rfloor_N }^N \bigr\|^2_{ L^p( \Omega; \mathbb{R}^d ) } \,ds
\end{eqnarray*}
for all $ t \in[0,T] $, $ N \in\mathbb{N}$ and all $ p \in[4, \infty)
$. In the next step Gronwall's lemma shows that
\begin{eqnarray*}
&& \Bigl\| {\sup_{ t \in[0,T] }} \| X_t - \bar{Y}_t^N \| \Bigr\|^2_{ L^p( \Omega;
\mathbb{R} ) }
\\[-2pt]
&&\qquad\leq
2 e^{ 4 T ( p^2 c^2 m + c + 1 ) } \biggl( \int_0^T \bigl\| \mu( \bar{Y}_s^N ) -
\mu\bigl( \bar{Y}_{ \lfloor s \rfloor_N }^N \bigr) \bigr\|^2_{ L^p( \Omega;
\mathbb{R}^d ) } \,ds
\\[-2pt]
&&\qquad\quad\hspace*{77.5pt}{} +
\frac{ T^2 }{ N^2 } \int_0^T \bigl\| \mu\bigl( \bar{Y}_{ \lfloor s \rfloor_N }^N
\bigr) \bigr\|^4_{ L^{2p}( \Omega; \mathbb{R}^d ) } \,ds
\\[-2pt]
&&\qquad\quad\hspace*{77.5pt}{} +
2 p^2 c^2 m \int_0^T \bigl\| \bar{Y}_s^N - \bar{Y}_{ \lfloor s \rfloor_N }^N
\bigr\|^2_{ L^p( \Omega; \mathbb{R}^d ) } \,ds \biggr)\vadjust{\goodbreak}
\end{eqnarray*}
and hence, the inequality $\sqrt{a+b+c}\leq\sqrt{a}+\sqrt{b}+\sqrt{c}$
for all $a,b,c\in[0,\infty)$ gives that
%
%
\begin{eqnarray}\label{eqineq1}\quad
&&\Bigl\| {\sup_{ t \in[0,T] }} \| X_t - \bar{Y}_t^N \| \Bigr\|_{ L^p(
\Omega; \mathbb{R} ) }
\nonumber\\
&&\qquad\leq
\sqrt{ 2T } e^{ 2 T ( p^2 c^2 m + c + 1 ) } \biggl( \sup_{ t
\in[0,T] } \bigl\| \mu( \bar{Y}_t^N ) - \mu\bigl( \bar{Y}_{ \lfloor t \rfloor_N
}^N \bigr) \bigr\|_{ L^p( \Omega; \mathbb{R}^d ) }
\nonumber\\[-8pt]\\[-8pt]
&&\qquad\quad\hspace*{94pt}{} +
\frac{ T }{ N } \Bigl[ {\sup_{ n \in\{ 0, 1, \ldots, N \} }}
\| \mu( Y^N_n ) \|^2_{ L^{2p}( \Omega; \mathbb{R}^d ) } \Bigr]
\nonumber\\
&&\qquad\quad\hspace*{94pt}{} +
p c \sqrt{2m} \Bigl[ \sup_{ t \in[0,T] } \bigl\| \bar{Y}_t^N -
\bar{Y}_{ \lfloor t \rfloor_N }^N \bigr\|_{ L^p( \Omega; \mathbb{R}^d ) } \Bigr]
\biggr)\nonumber
\end{eqnarray}
for all $ N \in\mathbb{N} $ and all $ p \in[4, \infty) $. Additionally,
the Burkholder--Davis--Gundy type inequality in Lemma~\ref{lBDG} shows
that
\begin{eqnarray*}
&&
\sup_{ t \in[0,T] } \bigl\| \bar{Y}_t^N - \bar{Y}_{ \lfloor t \rfloor_N }^N
\bigr\|_{ L^{p}( \Omega; \mathbb{R}^d ) }
\\
&&\qquad\leq\frac{T}{N} \biggl( \sup_{ t \in[0,T] } \biggl\| \frac{ \mu( \bar{Y}^N_{
\lfloor t \rfloor_N } ) }{ 1 + ({ T }/{N}) \| \mu( \bar{Y}^N_{
\lfloor t \rfloor_N } ) \| }\biggr\|_{ L^p( \Omega; \mathbb{R}^d ) } \biggr)
\\
&&\qquad\quad{}
+ \sup_{ t \in[0,T] } \biggl\| \int_{ \lfloor t \rfloor_N }^t \sigma\bigl(
\bar{Y}^N_{ \lfloor t \rfloor_N } \bigr) \,dW_s \biggr\|_{ L^p( \Omega;
\mathbb{R}^d ) }
\\
&&\qquad\leq
\frac{T}{ \sqrt{N} } \Bigl( \sup_{ n \in\{ 0, 1, \ldots, N \} } \| \mu(
Y_n^N ) \|_{ L^{p}( \Omega; \mathbb{R}^d ) } \Bigr)
\\
&&\qquad\quad{}
+ \frac{ p \sqrt{T m} }{ \sqrt{N} } \Bigl( {\sup_{i\in\{1,2,\ldots,m\}}
\sup_{ n \in\{ 0, 1, \ldots, N \} } }\| \sigma_i( Y_n^N ) \|_{ L^{p}(
\Omega; \R^{d} ) } \Bigr)
\end{eqnarray*}
for all $ N \in\mathbb{N} $ and all $ p \in[2, \infty) $. Lemma
\ref{lemmubound} hence implies that
%
%
\begin{equation}
\label{eqytdif} \sup_{ N \in\mathbb{N} } \Bigl( \sqrt{N} \Bigl[ \sup_{ t \in[0,T]
} \bigl\| \bar{Y}_t^N - \bar{Y}_{ \lfloor t \rfloor_N }^N \bigr\|_{ L^{p}(
\Omega; \mathbb{R}^d ) } \Bigr] \Bigr) < \infty
\end{equation}
for all $ p \in[1,\infty) $. In particular, we obtain that
%
%
\begin{equation}
\label{eqytbound} {\sup_{ N \in\mathbb{N} } \sup_{ t \in[0,T] }} \|
\bar{Y}_t^N \|_{ L^{p}( \Omega; \mathbb{R}^d ) } < \infty
\end{equation}
for all $ p \in[1,\infty) $ due to Lemma~\ref{cmomEuler}. Moreover, the
estimate $ \| \mu( x ) - \mu( y ) \| \leq c( 1 + \| x \|^c + \| y \|^c
) \| x - y \| $ for all $ x, y \in\mathbb{R}^d $ gives that
%
%
\begin{eqnarray}
\label{eqineq3}
&&
\sup_{ t \in[0,T] } \bigl\| \mu( \bar{Y}_t^N ) - \mu\bigl(
\bar{Y}_{ \lfloor t \rfloor_N }^N \bigr) \bigr\|_{ L^{p}( \Omega; \mathbb{R}^d )
}
\nonumber\hspace*{-35pt}\\[-8pt]\\[-8pt]
&&\qquad\leq c\Bigl( 1 + {2\sup_{ t \in[0,T] }} \| \bar{Y}_t^N \|^c_{ L^{2pc}(
\Omega; \mathbb{R}^d ) } \Bigr) \Bigl( \sup_{ t \in[0,T] } \bigl\| \bar{Y}_t^N -
\bar{Y}_{ \lfloor t \rfloor_N }^N \bigr\|_{ L^{2p}( \Omega; \mathbb{R}^d )
}\Bigr)\nonumber\hspace*{-35pt}
\end{eqnarray}
for all $ N \in\mathbb{N} $ and all $ p \in[1, \infty) $ and
inequalities~(\ref{eqytdif}) and~(\ref{eqytbound}) hence show that
%
%
\begin{equation}
\label{eqmuytdif} \sup_{ N \in\mathbb{N} } \Bigl( \sqrt{N} \Bigl[ \sup_{ t
\in[0,T] } \bigl\| \mu( \bar{Y}_t^N ) - \mu\bigl( \bar{Y}_{ \lfloor t \rfloor_N
}^N \bigr) \bigr\|_{ L^{p}( \Omega; \mathbb{R}^d ) } \Bigr] \Bigr) < \infty
\end{equation}
for all $ p \in[1,\infty) $. Combining~(\ref{eqineq1}),
(\ref{eqytdif}),~(\ref{eqmuytdif}) and Lemma~\ref{lemmubound}
finally\break
shows~(\ref{eqmainresult}). This completes the proof of Theorem
\ref{thmmain}.

\section*{Acknowledgments}

The authors thank an anonymous referee for very helpful
comments.


%
%

\printaddresses

\end{document}